\documentclass[letterpaper, 10 pt, conference]{ieeeconf}  % Comment this line out if you need a4paper

\IEEEoverridecommandlockouts                              % This command is only needed if 
\usepackage{amsmath}
\usepackage{bm}        % for bold math symbol
\usepackage{amssymb}
\usepackage{color}

\usepackage{graphicx}      % include this line if your document contains figures
\usepackage{tikz}
\usepackage{makecell} % for newline in a cell
\usepackage{stfloats} % to make 2 column tables stay in the bottom of the page
\usepackage{bbold}
\usepackage{cuted}
\usepackage[T1]{fontenc} % for \guillemotleft
\usepackage{wasysym} % for \Square
\usepackage{tabularx}
\usepackage{multirow}
\usepackage{accents}
\usepackage{comment}
\usepackage{mathtools}
\usepackage{caption}
\usepackage{subcaption}

\newtheorem{assumption}{Assumption}

\title{\LARGE \bf
Trajectory Optimization and NMPC Tracking for a Fixed--Wing UAV in Deep Stall with Perch Landing*
}

\author{Huu Thien Nguyen$^{1}$, Ionela Prodan$^{2}$, and Fernando A. C. C. Fontes$^{1}$% <-this % stops a space
\thanks{*The authors thank the valuable contributions of Prof. Fernando Lobo Pereira during several stages of this research. The work of Huu Thien Nguyen is funded by Fundação para a Ciência e a Tecnologia (FCT) Portugal, under the contract number $2020.07959.BD$.}% <-this % stops a space
\thanks{$^{1}$Huu Thien Nguyen and Fernando A. C. C. Fontes are with SYSTEC -- ISR, University of Porto, Porto, Portugal.
        {\tt\small \{nguyen,\ faf\}@fe.up.pt}.}%
\thanks{$^{2}$Ionela Prodan is with Univ. Grenoble Alpes, Grenoble INP, LCIS, F-26000, Valence, France.
        {\tt\small ionela.prodan@lcis.grenoble-inp.fr}.}%
}

\begin{document}

\maketitle
\thispagestyle{empty}
\pagestyle{empty}

%%%%%%%%%%%%%%%%%%%%%%%%%%%%%%%%%%%%%%%%%%%%%%%%%%%%%%%%%%%%%%%%%%%%%%%%%%%%%%%%
\begin{abstract}
This paper presents a novel recovery technique for a fixed-wing UAV (Unmanned Aerial Vehicle) based on constrained optimization: i) we propose a trajectory generation for landing the UAV where it first reduces its altitude by deep stalling, then perches on a recovery net, ii) we design an NMPC (Nonlinear Model Predictive Control) tracking controller with terminal constraints for the optimal generated trajectory under disturbances. Compared to nominal net recovery procedures, this technique greatly reduces the landing time and the final airspeed of the UAV.
Simulation results for various wind conditions demonstrate the feasibility of the idea. 
\end{abstract}

%%%%%%%%%%%%%%%%%%%%%%%%%%%%%%%%%%%%%%%%%%%%%%%%%%%%%%%%%%%%%%%%%%%%%%%%%%%%%%%%
\section{INTRODUCTION}

Landing fixed-wing (FW) UAVs is a challenging task. Unlike multicopters with propellers intentionally positioned for a safe vertical landing, many FW UAVs do not have a built-in mechanism (e.g., a landing gear) to dampen the impact when touching the ground. Hence, the deceleration at the moment of impact negatively affects the mechanical structure of the UAVs. For this kind of UAVs, parachutes, nets, and wires are often used as recovery techniques \cite{skitmore_LaunchRecovery_2020}. However, these auxiliary devices entail operational constraints such as the need for a priori deployment of infrastructures, and restrictions on the landing area \cite{ klausen_AutonomousRecovery_2018}. It is worth noting that, in the aforementioned reference, two multicopters are used to hang the net and move along with the UAV to reduce the impact force. 

In typical operations, FW UAVs must stay outside the stall region, where the Angle of Attack (AoA) provides the largest lift. Above this value, defined as the ``critical AoA'', the aircraft falls into the ``post-stall'' mode of operation, in which the lift is lost, the controllability is reduced, and the drag is increased. However, by adopting appropriate control strategies, significant operational advantages can be extracted by making use of this large AoA region. Thus, the existence of several research works addressing these challenges is not surprising. 

Deep stall happens when the aircraft surpasses its critical AoA, the airflow surrounding the wings separates. When the airflow returns to stable, the aircraft dives in the post-stall region \cite{taniguchi_AnalysisDeepstall_2008}. Perching, on the other hand, is a technique inspired by nature which also exploits the separation airflow, and the high drag force but at a higher AoA $(>90^\circ)$ to land an aircraft at a sufficiently small airspeed \cite{wickenheiser_OptimizationPerching_2008,alikhan_FlightDynamics_2013}.

Previous works on deep stall include \cite{kim_VisionAssisted_2022}, where the UAV is vision guided with the help of the computer mouse and two PI controllers are used to land the UAV.
% In \cite{park_ControlGuidance_2020}, trim analysis in the deep stall phase is performed by linearizing around trim points, and guidance control laws are designed to land headwind. 
Cunis et al. \cite{cunis_DynamicStability_2020} use bifurcation analysis to analyze the dynamic stability of a UAV in the deep stall and post-stall region.
Mathisen et al. \cite{mathisen_NonlinearModel_2016} propose a deep stall landing procedure using an NMPC controller to land an FW UAV on a predefined point in three-dimensional space. Extensive simulations are executed to demonstrate the relation between the wind velocity and the flight path angle. 
The algorithm is augmented in \cite{mathisen_PrecisionDeepStall_2021} with software-in-the-loop simulations. 

Regarding works on perch, \cite{venkateswararao_ParametricStudy_2015} proposes an optimization problem to minimize the distance traveled while constraining the final airspeed to be less than $5\%$ of the initial one but only with the longitudinal dynamics.
Feroskhan et al. \cite{feroskhan_SolutionsPlanar_2020} follow this cost formulation and final airspeed constraint to perch a UAV in three dimensions.
Reinforcement learning is used in \cite{waldock_LearningPerform_2018} to generate perching trajectories for a variable-sweep wing UAV.
% In \cite{wickenheiser_OptimizationPerching_2008}, an optimization problem with two cost functions is solved: minimization of the starting distance, and of the altitude overshoot for perching by controlling a morphing UAV. 
% In \cite{alikhan_FlightDynamics_2013}, an optimal control problem (OCP) minimizing the total distance while constraining the final states to get the perching landing trajectory for a fixed-wing aircraft in three dimensions is considered.
Moore et al. \cite{moore_RobustPoststall_2014} use LQR-Trees algorithm to robustly perch an FW glider. 
The authors in \cite{cory_ExperimentsFixedWing_2008}, from experiments, find out that the flat-plate model is well-suited for the operation in the post-stall regime.

To the best of our knowledge, no work on UAV deep stall landing has carefully considered the final airspeed of the UAV at the moment before touching the ground, and its effect on landing performance. Thus, in this article, we address this issue, and design a landing strategy on a recovery net, whose performance compares with the current alternatives as shown in Table \ref{table_comparison_landing_techniques} (more $+$ signs means a larger value). A previous version of this article can be found online \cite{nguyen_TrajectoryOptimization_2022}.

\begin{table}[ht!]
	\caption{Comparison of landing techniques for FW UAVs}
	\label{table_comparison_landing_techniques}
	\begin{center}
		\begin{tabular}{cccc}
			\textbf{Technique}           & \textbf{Altitude change} & \textbf{Final airspeed} & \textbf{Landing time}  \\
			\hline
			Deep stall   & $+++$               & $+++$            & $+$      \\
			\hline
			Perch        & $+$                 & $+$              & $+$    \\
			\hline
	Net recovery & $+$                 & $++$              & $++$    \\
			\hline
			Our approach & $+++$              & $++$             & $++$   
		\end{tabular}    
	\end{center}
\end{table}

The contributions in this paper are summarized as follows:
\begin{itemize}
	\item It combines deep stall with perch landing to achieve both short landing time in a narrow space while preserving a small final landing airspeed. 
	\item A deep stall with perch landing reference trajectory is generated  by solving a constrained Optimal Control Problem (OCP). 
	\item An NMPC tracking controller is developed, allowing for deep stall landing and perching of the UAV under windy conditions.
	\item \textcolor{black}{The feasibility of the overall scheme is shown by simulations under several wind conditions.}
\end{itemize}

Notation: For an arbitrary vector $\mathbf{x}$, $\|\mathbf{x}\|^2_P = \mathbf{x}^\intercal P \mathbf{x}$. Let $\mathbb{I}_n$ represent the identity matrix of size $n$, $\mathbb{S}^n_{+}$ ($\mathbb{S}^n_{++}$) denote the vector space of $n\times n$ real symmetric positive semidefinite (positive definite) matrices. The subscript $\square_r$ denotes the reference values, while $\bar{\ast}$ and $\underbar{$\ast$}$ are the upper-bound and lower-bound defined for the variable $\ast$, respectively.

This article is organized as follows: in Section \ref{sec_Longitudinal_dynamics}, the longitudinal dynamics of an FW UAV is presented. The problem formulation in Section \ref{sec_Problem_Formulation} includes the OCP to generate the reference trajectory, and the NMPC scheme to track it. Simulation results  based on the specific data of the Aerosonde UAV are presented in Section \ref{sec_Simulation}. Finally, some brief conclusions and prospective future work are outlined in Section \ref{sec_Conclusions}.
%%%%%%%%%%%%%%%%%%%%%%%%%%%%%%%%%%%%%%%%%%%%%%%%%%%%%%%%	
\section{Longitudinal fixed-wing UAV dynamics}
\label{sec_Longitudinal_dynamics}
\begin{figure}[ht!]
	\begin{center}		 	
		\includegraphics[width=0.4\columnwidth]{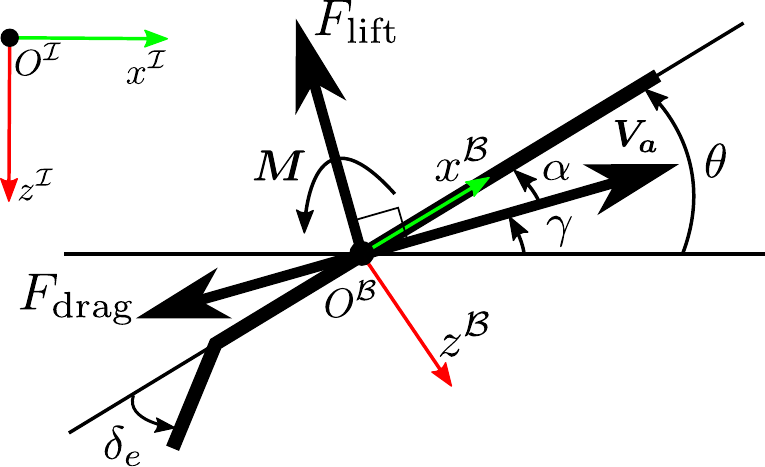}	
		\caption{Aerodynamic forces and moment acting on a longitudinal fixed-wing UAV}
		\label{fig_Aerosonde_Aerodynamic_Forces_Moment}
	\end{center}
\end{figure}

%%%%%%%%%%%%%%%%%%%%%%%%%%%%%%%%%%%%%%%%%%%%%%%%%%%%%%%%
We consider two reference frames: the inertial frame $\mathcal{I}\{O^\mathcal{I},x^\mathcal{I},z^\mathcal{I}\}$ fixed, pointing north and down; the body frame $\mathcal{B}\{O^\mathcal{B},x^\mathcal{B},z^\mathcal{B}\}$ attached to the center of mass of the UAV which moves along with the UAV, as in Fig. \ref{fig_Aerosonde_Aerodynamic_Forces_Moment}.
The two-dimensional longitudinal dynamics of an FW UAV \cite{beard_SmallUnmanned_2012} are as follows 
\begin{equation}
		\label{eqn_uav_2d_dynamics}
\begin{aligned}
	\begin{bmatrix}
		\dot{x}\\
		\dot{z}
	\end{bmatrix}
	&=
	\mathcal{R}^{\mathcal{I}}_{\mathcal{B}}(\theta)
	\begin{bmatrix}
		u\\
		w
	\end{bmatrix} \\
	\dot{u} &=-\omega_y w + {f}_{x}/m \\
	\dot{w} &= \omega_y u +  {f}_{z}/m \\
	\dot{\theta} &= \omega_y \\
	\dot{\omega}_y &=	{M}/{J_{y}}
\end{aligned}
\end{equation}
\begin{comment}
\begin{equation}
	\label{eqn_uav_2d_dynamics}
	\begin{bmatrix}
		\dot{x} \\
		\dot{z} \\
		\dot{u} \\
		\dot{w} \\
		\dot{\theta} \\
		\dot{\omega}_y
	\end{bmatrix}
	=
	\begin{bmatrix}
		\cos (\theta) u +\sin (\theta) w \\
		-\sin (\theta) u +\cos (\theta) w \\
		-\omega_y w + {f}_{x}/m \\
		\omega_y u +  {f}_{z}/m \\
		\omega_y \\
		{M}/{J_{y}}
	\end{bmatrix},
\end{equation}
\end{comment}
where $x$ and $z$ are the horizontal and vertical position of the UAV in $\mathcal{I}$; $\mathcal{R}^{\mathcal{I}}_{\mathcal{B}}(\theta) \triangleq \begin{bsmallmatrix}
	\cos (\theta) & \sin (\theta) \\
	-\sin (\theta) & \cos (\theta)
\end{bsmallmatrix}$ is the rotation matrix from $\mathcal{B}$ to $\mathcal{I}$, satisfying $\mathcal{R}^{\mathcal{I}}_{\mathcal{B}}(\theta) =[\mathcal{R}^{\mathcal{B}}_{\mathcal{I}}(\theta)]^{-1}$; $u$ and $w$ are longitudinal and vertical velocity of $\mathcal{B}$ w.r.t $\mathcal{I}$, expressed in $\mathcal{B}$; $\theta$ is the pitch angle; $\omega_y$ is the pitch rate in $\mathcal{B}$; $m$ is the mass of the UAV; $J_y$ is the moment of inertia about $y^\mathcal{B}$; ${f}_{x}$ and ${f}_{z}$ are the externally applied forces in $\mathcal{B}$; ${M}$ is the pitching moment around $y^\mathcal{B}$. 

The aerodynamic forces and moment acting on the UAV, also shown in Fig. \ref{fig_Aerosonde_Aerodynamic_Forces_Moment}, are given by 
\begin{subequations}
	\label{eqn_lift_drag_moment}
	\begin{align}
		\label{eqn_f_lift_drag}
	&	{F}_i =  \frac{1}{2} \rho V_{a}^2 S \left[ C_{i}(\alpha)+C_{i_{q}} \frac{c}{2 V_{a}} \omega_y + C_{i_{\delta_{e}}} \delta_{e} \right], \\
		\label{eqn_moment}
	 &	M= \frac{1}{2} \rho V_{a}^2 S c\left[C_m(\alpha) + C_{m_{q}} \frac{c}{2 V_{a}} \omega_y +C_{m_{\delta_{e}}} \delta_{e}\right],
	\end{align}
\end{subequations}
where $i \in \{L,D\}$ represents the lift and drag, $\rho$ is the air density, $S$ is the planform area, $c$ is the mean chord of the wing, $C_{i_{q}}$, $C_{i_{\delta_{e}}}$, $C_{m_{q}}$, $C_{m_{\delta_{e}}}$ are aerodynamic constants, $\delta_e$ is the elevator deflection angle, while $V_a$ and $\alpha$ are the airspeed and the angle of attack (AoA) calculated as follows
\begin{subequations}
	\label{eqn_airspeed_and_AoA}
	\begin{align}
		\label{eqn_airspeed}
		V_{a} =\sqrt{u_{a}^{2}+w_{a}^{2}}, \ \
		& \alpha = \tan^{-1} \left( {w_a}/{u_a} \right), \\
		\label{eqn_state_3_state_4_wind}
		\text{with }\begin{bmatrix}
			u_{a} \\
			w_{a}
		\end{bmatrix} 
		= &
		\begin{bmatrix}
			u - u_{w} \\
			w - w_{w}
		\end{bmatrix}.
	\end{align}
\end{subequations}
Here, $u_{w}$ and $w_{w}$ are the wind components in $\mathcal{B}$, whose positive values mean the vector components of the wind have the same direction with the $x$ and $z$ axes in $\mathcal{B}$. Furthermore, $C_i(\alpha)$ $(i \in \{L,D\})$, $C_m(\alpha)$ in \eqref{eqn_lift_drag_moment} are lift, drag, and pitching moment coefficients, respectively
\begin{subequations}
	\label{eqn_C_L_C_D_C_m}
	\begin{align}
	&	C_{i}(\alpha) = [1-\sigma(\alpha)] C_{i}^{pre}(\alpha) + \sigma(\alpha) C_{i}^{post}(\alpha) , \quad \\
	&	C_m(\alpha) = C_{m_0} +  C_{m_\alpha} \alpha.
	\end{align}
\end{subequations}
When the UAV is in the pre-stall regime, the typical linear and quadratic functions for lift and drag are used \cite{beard_SmallUnmanned_2012}, while the post-stall aerodynamics are inherited from the flat-plate model in \cite{cory_ExperimentsFixedWing_2008}
\begin{subequations}
	\begin{align}
		C_{L}^{pre}(\alpha) &= C_{L_0} + C_{L_\alpha} \alpha ,\qquad C_{L}^{post}(\alpha) = \sin(2 \alpha),  \\
		C_{D}^{pre}(\alpha)&=C_{D_p} + \tfrac{[	C_{L}^{pre}(\alpha)]^2}{\pi e AR}  ,\ C_{D}^{post}(\alpha)= 2 \sin^2(\alpha).  
	\end{align}
\end{subequations}
Here, $C_{D_p}$ is the parasitic drag, $e$ is the Oswald efficiency factor, $AR$ is the wing aspect ratio, the AoA $\alpha$ has the unit $[rad]$, and the $\sigma(\alpha)$ in (4a) is a sigmoid function used as a blending function between the two regimes
\begin{equation}
	\sigma(\alpha) = \frac{1}{1+e^{-\tilde{M}(\alpha-\alpha_0)}}
\end{equation}
where $\tilde{M}$ is the transition rate and $\alpha_0$ is the cut off AoA.
%From Fig. \ref{fig_Aerosonde_Aerodynamic_Forces_Moment}, the pitch angle $\theta$ is the angle between the aircraft's longitudinal axis w.r.t. the horizontal line, but the direction that the UAV pointing to is actually ruled by the flight path angle $\gamma$, satisfying the relation
From Fig. \ref{fig_Aerosonde_Aerodynamic_Forces_Moment}, the pitch angle $\theta$, the flight path angle $\gamma$, and the AoA $\alpha$ satisfying the relation
\begin{equation}
	\label{eqn_pitch_AoA_gamma}
	\theta = \alpha + \gamma.
\end{equation}
Moreover, the external forces ${f}_{x}$ and ${f}_{z}$ in \eqref{eqn_uav_2d_dynamics} consist of
\begin{subequations}
	\begin{align}
		\begin{bmatrix}
			{f}_{x},
			{f}_{z}
		\end{bmatrix}^\intercal &= \bm{f_g}+\bm{f_a}+\bm{f_p}, \\
		\bm{f_g} &=
		\begin{bmatrix}
			- \sin \theta,
			\cos \theta
		\end{bmatrix}^\intercal
		mg, \\
		\label{eqn_forces_xz_aerodynamics}
		\bm{f_a} &=
		\begin{bmatrix}
			\sin{\alpha} & -\cos{\alpha} \\
			-\cos{\alpha} & -\sin{\alpha}
		\end{bmatrix}
		\begin{bmatrix}
			{F}_L \\
			{F}_D
		\end{bmatrix}, \\
		\bm{f_p} &=\frac{1}{2} \rho S_{\text {prop }} C_{\text {prop }}
		\begin{bmatrix}
			\left(k_{\text {motor }} \delta_{t}\right)^{2}-V_{a}^{2},
			0
		\end{bmatrix}^\intercal \label{eqn_f_p},
	\end{align}
\end{subequations}
where $\bm{f_g}$ is the gravitational force, $\bm{f_a}$ is the aerodynamic force in \eqref{eqn_f_lift_drag}, $\bm{f_p}$ is the propulsive force, $S_{\text {prop }}$ and $C_{\text {prop }}$ are propeller's parameters, $k_{\text {motor }}$ is the motor constant, and $\delta_t \in [0,1]$ is the pulse-width modulation of the propeller.

Therefore, the UAV dynamics \eqref{eqn_uav_2d_dynamics} can be written in the canonical nonlinear form as follows
\begin{equation}
	\dot{\mathbf{x}} = f(\mathbf{x},\mathbf{u}),
\end{equation}
where $\mathbf{x} \triangleq \left[ x,z,u,w,\theta,\omega_y \right]^\intercal$, and $\mathbf{u} \triangleq \left[ \delta_{e}, \delta_{t} \right]^\intercal$ gather the states and the inputs of the system, respectively.

%%%%%%%%%%%%%%%%%%%%%%%%%%%%%%%%%%%%%%%%%%%%%%%%%%%%%%%%%%%%%%%%%%%%%%%%%%%
%%%%%%%%%%%%%%%%%%%%%%%%%%%%%%%%%%%%%%%%%%%%%%%%%%%%%%%%%%%%%%
\section{PROBLEM FORMULATION}
\label{sec_Problem_Formulation}
The problem is stated as follows. The procedure to autonomously deep stall with perch land an FW UAV consists of two tasks. While a reference trajectory and the associated inputs are generated in the first task, the controller for the UAV to track the generated trajectory is designed in the second task. For the sake of simplicity, the deep stall and the perch maneuvers are combined in one landing phase. Let $t_0$ and $t_f$ be the initial time and the final time when the UAV touches a recovery net in the landing phase. The recovery net is fixed at the horizontal position $x_\text{net}=0$, with its height bounded in $[ \underaccent{\bar}{z}_\text{net} ,\bar{z}_\text{net}]$, \textcolor{black}{where $\underaccent{\bar}{z}_\text{net}$, $\bar{z}_\text{net}$ are negative, and }the magnitude of $\bar{z}_\text{net}$ equals to the length of the UAV \textcolor{black}{(see Fig. \ref{fig_LOS_pos})}. The position where the UAV initiates its landing maneuver has the coordinates $(x_0,z_0)$. 

\begin{assumption}
Before entering the landing phase, the UAV cruises in steady level flight, that is, for the existing wind condition, it satisfies\footnote{The conditions $\theta = \alpha$ and $\dot{z}=0$ are due to $\gamma=0^\circ$.}
\begin{equation}
	\label{eqn_level_trim_flight}
	\dot{u}=0, \dot{w}=0, \theta = \alpha, \dot{z}=0, \omega_y=0, V_a = V_{a_r}.
\end{equation}
\end{assumption}
% to obtain the set of trim states and inputs denoted with the subscript $\square_0$: $u_0$, $w_0$, $V_{a_0}=V_{a_r}$, $\theta_0=\alpha_0$, $\omega_{y_0}=0$, $ \delta_{e_0}$, and $ \delta_{t_0}$. 
Together with the initial position coordinate, we have the initial trim state and input, denoted with the subscript $\square_0$, used hereinafter for the trajectory generation phase
\begin{subequations}
	\label{eqn_initial_state_input_trim}
	\begin{align}
\mathbf{x}_0 & \triangleq [x_0,z_0,u_0,w_0,\theta_0, \omega_{y_0}]^\intercal, \\
\mathbf{u}_0 & \triangleq [\delta_{e_0}, \delta_{t_0}]^\intercal.
	\end{align}
\end{subequations}

% \textcolor{black}{
% \begin{assumption}
% There always exists a solution for \eqref{eqn_level_trim_flight} satisfying the physical constraints of the FW UAV.
% \end{assumption}
% }
%%%%%%%%%%%%%%%%%%%%%%%%%%%%%%%%%%%%%%%%%%%%%%%%%%%%%%%%%%%%%%
\subsection{Control requirement for the net recovery}
\begin{figure}[ht!]
		\centering
		\includegraphics[width=0.9\columnwidth]{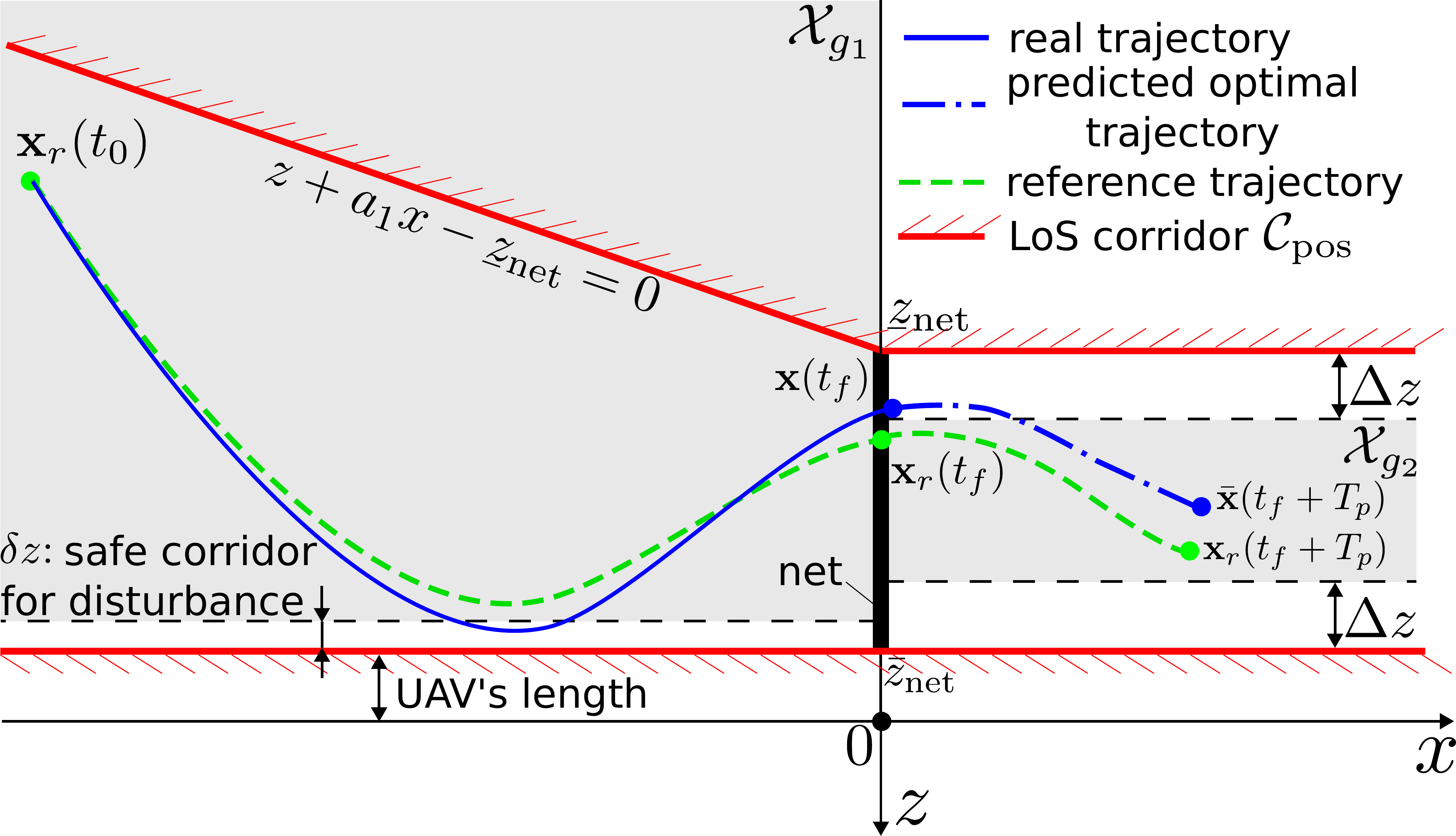}
		\caption{$\mathcal{X}_{g_1}$, $\mathcal{X}_{g_2}$, and $\mathcal{C}_{\text{pos}}$ projected to the $xz$ plane.}
		\label{fig_LOS_pos}
\end{figure}

\begin{figure}[ht!]
		\centering
		\includegraphics[width=0.9\columnwidth]{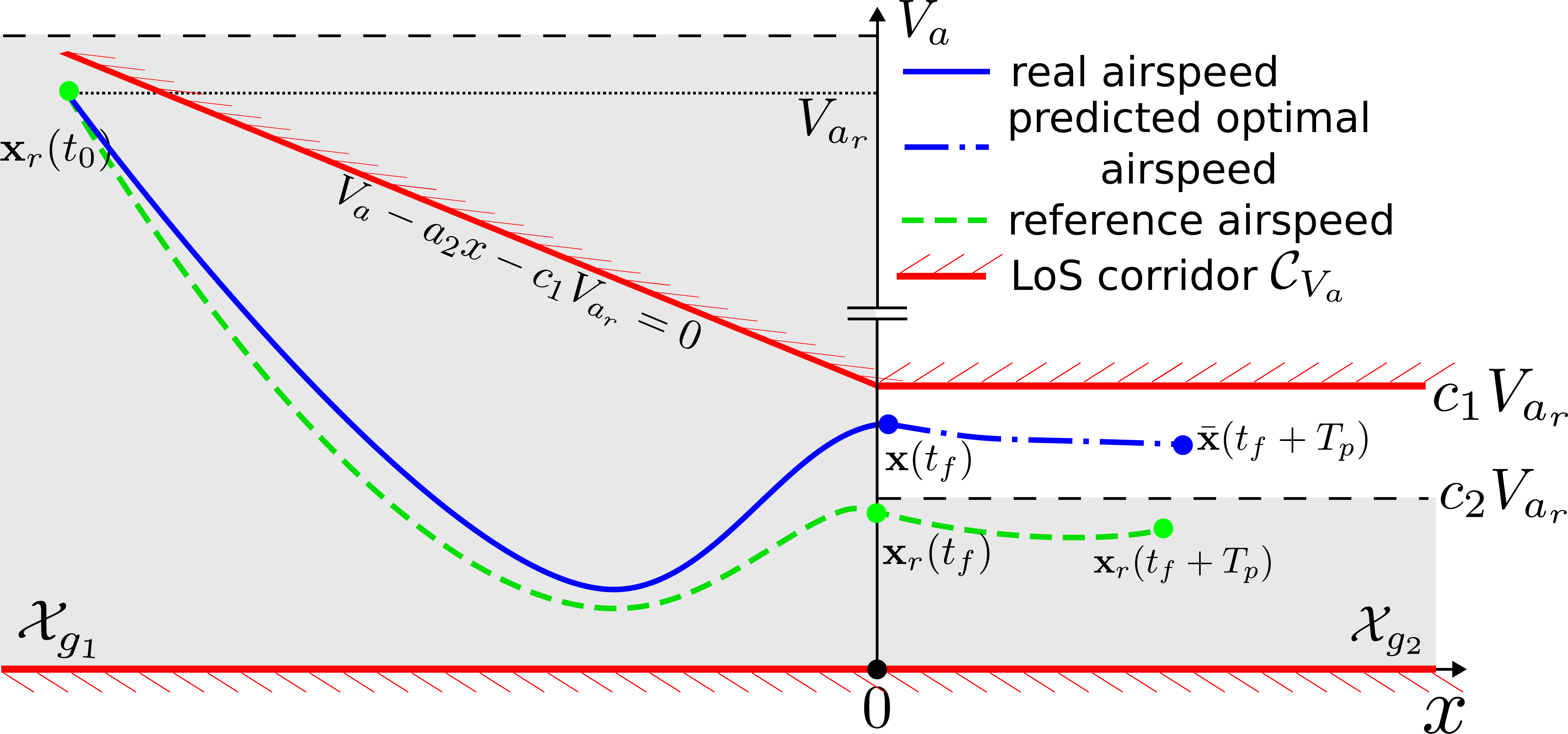}
		\caption{$\mathcal{X}_{g_1}$, $\mathcal{X}_{g_2}$, and $\mathcal{C}_{V_a}$ projected to the $x V_a$ plane.}
		\label{fig_LOS_airspeed}
\end{figure}

We define that a successful landing on the net at time $t_f$ near horizontal position $x=0$, with sufficiently small airspeed, is the one satisfying the following requirements
\begin{subequations}
	\label{eqn_requirements}
	\begin{align}
		V_{a}(t_f) & \leq c_1  V_{a_r}, \label{eqn_req_track_Va} \\
				\underaccent{\bar}{z}_\text{net}   \leq z_r(t_f) & \leq \bar{z}_\text{net}, \label{eqn_req_track_zf} \\
		-\dot{x}(t_f) \delta t \leq x(t_f) & \leq \dot{x}(t_f) \delta t. \label{eqn_req_track_xf}
	\end{align}
\end{subequations}
The condition \eqref{eqn_req_track_Va} is to classify the landing as perching, where $c_1$ is a positive constant. The condition \eqref{eqn_req_track_zf} is to verify the UAV falls into the recovery net. Furthermore, the condition \eqref{eqn_req_track_xf} is the horizontal tolerance, $\delta t$ in \eqref{eqn_req_track_xf} is the sampling time that will be explained in Section \ref{subsec_Landing_Traject_Track}. These constraints need to be satisfied in real life situations, where various uncertainty such as wind gusts may manifest. Therefore, we define a target reference trajectory that complies with stricter requirements
\begin{subequations}
	\begin{align}
		V_{a_r}(t_f) & \leq c_2  V_{a_r}, \label{eqn_req_gen_Va} \\
				\underaccent{\bar}{z}_\text{net} + \Delta z  \leq z_r(t_f) & \leq \bar{z}_\text{net} - \Delta z, \label{eqn_req_gen_zf} \\
		 0 \leq x_r(t_f) & \leq \delta x, \label{eqn_req_gen_xf}
	\end{align}
\end{subequations}
where $0<c_2<c_1$, $ \Delta z>0$, and \textcolor{black}{ $\delta x$ is a small positive number}. These requirements will be imposed as path-wise and terminal state constraints in the optimal control problems by inclusion in the sets $\mathcal{X}_{g_1}$, $\mathcal{X}_{g_2}$, and $\mathcal{X}_{f} $ as in \eqref{eqn_set_X_g_1}--\eqref{eqn_set_X_f} introduced below. We define the following sets 
\begin{subequations}
	\label{eqn_set_constraints}
	\begin{align}
		\mathcal{U} = &
		\{ \mathbf{u} \in \mathbb{R}^2  |  
		\left[ 
		\underaccent{\bar}{\delta}_{e} , 
		\underaccent{\bar}{\delta}_{t} \right]^\intercal 
		\leq 
		\left[ {\delta_{e} , \delta_{t}} \right]^\intercal 
		\leq 
		\left[\bar{\delta}_e , 	\bar{\delta}_t \right]^\intercal \}, 
		\label{eqn_set_U} \\
		\mathcal{X}_{g_1}  = &
		\{ \mathbf{x} \in \mathbb{R}^6 | 
		x \leq 0, \allowbreak 
		 z \leq \bar{z}_\text{net} \hspace{-0.1cm} - \hspace{-0.1cm} \delta z, 
		 \mathrm{atan}\begin{pmatrix}\tfrac{w-w_w}{u-u_w}\end{pmatrix} \in [\underaccent{\bar}{\alpha},\bar{\alpha}],\nonumber \\
	& 	|u| \leq  \bar{u},  
		|w| \leq \bar{w},  
	  \theta \in [0,\bar{\theta}], 	
		|\omega_y| \leq \bar{\omega}_{y} \},  
		\label{eqn_set_X_g_1}\\
		\mathcal{X}_{g_2}  = &
		\{ \mathbf{x} \in \mathbb{R}^6 |  
		x \geq 0,  	z \in [	\underaccent{\bar}{z}_\text{net} + \Delta z ,  \bar{z}_\text{net} - \Delta z ],\nonumber \\	
		&	\sqrt{(u-u_w)^2+  (w-w_w)^2} \leq c_2 V_{a_r},  0^\circ \leq \theta \leq 90^\circ  ,  \nonumber \\
		&		90^\circ \leq \mathrm{atan}\begin{pmatrix}\tfrac{w-w_w}{u-u_w}\end{pmatrix} \leq \bar{\alpha},
		|\omega_y| \leq \bar{\omega}_{y} \},  
		\label{eqn_set_X_g_2}\\
\mathcal{X}_{f}=&  \mathcal{X}_{g_2} \cap \{\mathbf{x} \in \mathbb{R}^6 | x \in [0,\delta x] \},
\label{eqn_set_X_f}\\
	\mathcal{X}_\text{MPC}  = &
		\{ \mathbf{x} \in \mathbb{R}^6 | (x,z) \in \mathcal{C}_{\text{pos}}, \mathrm{atan}\begin{pmatrix}\tfrac{w-w_w}{u-u_w}\end{pmatrix} \in [\underaccent{\bar}{\alpha},\bar{\alpha}], \nonumber \\
		&		\theta \in [\underaccent{\bar}{\theta}_\text{MPC},\bar{\theta}_\text{MPC}], 
		|\omega_y| \leq \bar{\omega}_{y_\text{MPC}}, \nonumber \\
		 &	(x,\sqrt{(u-u_w)^2+  (w-w_w)^2}) \in \mathcal{C}_{V_a}  \}. 
		\label{eqn_set_X_MPC}
	\end{align}
\end{subequations} 
here, 
\textcolor{black}{$\delta_e$ and $\delta_t$ in \eqref{eqn_set_U} are the control inputs in \eqref{eqn_lift_drag_moment} and \eqref{eqn_f_p}; 
	$\delta z$ in \eqref{eqn_set_X_g_1} is a positive constant denoting the safe altitude when disturbances arise at the tracking task (see Fig. \ref{fig_LOS_pos} for $x \leq 0$);
	in \eqref{eqn_set_X_g_2}, $\Delta z$ is from \eqref{eqn_req_gen_zf}, $c_2$ is from \eqref{eqn_req_gen_Va}; 
	$\delta x$ in \eqref{eqn_set_X_f} is from \eqref{eqn_req_gen_xf};
	$u_w$ and $w_w$ as in \eqref{eqn_set_X_g_1}--\eqref{eqn_set_X_MPC} are from \eqref{eqn_state_3_state_4_wind};
}
$\mathcal{C}_{\text{pos}}$ and $\mathcal{C}_{V_a}$ in \eqref{eqn_set_X_MPC} are the \emph{Line-of-Sight (LoS) corridors} on $z$ and $V_a$ that we adopt the idea from \cite{dicairano_ModelPredictive_2012}
	\begin{equation}
		\mathcal{C}_{\text{pos}} = \left\{ x,z \middle|
		\begin{aligned}
			-z + \underaccent{\bar}{z}_\text{net} -a_1 x \leq & 0 \ \text{if} \ x<0\\
			-z +  \underaccent{\bar}{z}_\text{net} \leq & 0 \ \text{if} \ x \geq 0\\
			z - \bar{z}_\text{net} \leq  & 0 
		\end{aligned} \right\},
	\end{equation}
	\begin{equation}
		\mathcal{C}_{V_a} = \left\{ x,V_a \middle|
		\begin{aligned}
			V_a - c_1 V_{a_r} -a_2 x \leq & 0 \ \text{if} \ x < 0 \\
			V_a - c_1 V_{a_r} \leq &  0 \ \text{if} \ x \geq 0
		\end{aligned} \right\},
	\end{equation}
	where the slopes $a_1$ and $a_2$ defined after obtaining the reference trajectory. \textcolor{black}{The sets $\mathcal{X}_{g_1}$, $\mathcal{X}_{g_2}$ are the gray regions, and the LoS corridors $\mathcal{C}_{\text{pos}}$, $\mathcal{C}_{V_a}$, are the regions inside the red lines in Figures \ref{fig_LOS_pos} and \ref{fig_LOS_airspeed}}, in which \textcolor{black}{the obvious condition $V_a \geq 0 \ \forall x$ from \eqref{eqn_airspeed} is implicitly imposed in $\mathcal{C}_{V_a}$.} 
%%%%%%%%%%%%%%%%%%%%%%%%%%%%%%%%%%%%%%%%%%%%%%%%%%%%%%%%%%%%%%
%%%%%%%%%%%%%%%%%%%%%%%%%%%%%%%%%%%%%%%%%%%%%%%%%%%%%%%%%%%%%%
\subsection{Landing trajectory generation}
The reference landing phase is obtained by solving the OCP
\begin{equation}
	\label{eqn_OC_problem_formulation}
	\underset{{\mathbf{u}}(\cdot)}{\text{min}}     \ J_1\left(\mathbf{x}, \mathbf{u}\right)
	= \int_{t_0}^{t_f +T_p } \ell_1(({\mathbf{x}}(t), {\mathbf{u}}(t)) dt ,
\end{equation}
\vspace{-15pt}
\begin{subequations}
	\label{eqn_OC_problem_formulation_constraints}
	\begin{align}     
\text{subject to: }		& \dot{\mathbf{x}}(t) = f(\mathbf{x}(t),\mathbf{u}(t)), \  t \in [t_0,t_f+T_P] \\
		\label{eqn_OC_problem_formulation_state_input_constraints}
		&   \mathbf{u}(t) \in \mathcal{U}, \  t \in [t_0,t_f+T_P]  \\
		\label{eqn_OC_problem_formulation_initial_states}
		& \mathbf{x}(t) \in \mathcal{X}_{g_1},  \  t \in [t_0,t_f) \\
		& \mathbf{x}(t) \in \mathcal{X}_{g_2},  \  t \in [t_f,t_f+T_p]\\
		& \mathbf{x}(t_0) = \mathbf{x}_0, \ \mathbf{u}(t_0) = \mathbf{u}_0, \\
		& \mathbf{x}(t_f) = \mathbf{x}_f \in \mathcal{X}_f \subset \mathcal{X}_{g_2}.
	\end{align}
\end{subequations}
Here, $\ell_1(({\mathbf{x}}(t), {\mathbf{u}}(t))$ is the cost function defined as
\begin{equation}
	\label{eqn_OC_generation_cost}
	\ell_1(({\mathbf{x}}(t), {\mathbf{u}}(t)) =
	\|\delta_t\|^2_{P_{\delta_t}} + \|\Delta {\mathbf{u}}\|^2_{P_{\Delta U}} ,
\end{equation}
The prediction horizon $T_p$ of the MPC in the tracking task which will be explained in the next Subsection \ref{subsec_Landing_Traject_Track}. 
The reference trajectory in $[t_f,t_f+T_P]$ is called the \emph{augmented reference trajectory}, which is used to guide the terminal states of the optimal trajectories in the tracking task to stay inside the LoS corridors \textcolor{black}{(see Figs. \ref{fig_LOS_pos} and \ref{fig_LOS_airspeed} for $x \geq 0$)}.
Note that, $\mathbf{x}_0$ and $\mathbf{u}_0$ are the initial state and input of the landing phase as in \eqref{eqn_initial_state_input_trim}, which contains the solution of \eqref{eqn_level_trim_flight} for a predefined initial airspeed, $\mathbf{x}_f$ is the final landing state when the UAV perches on the recovery net. 
$P_{\delta_t} \in \mathbb{R}^+$ and ${P_{\Delta U}} \in \mathbb{S}^2_{++}$ are weighting scalars and matrix. $\mathcal{U}$, $\mathcal{X}_{g_1}$, $\mathcal{X}_{g_2}$, and $\mathcal{X}_{f}$ are as in \eqref{eqn_set_U}--\eqref{eqn_set_X_f}.

Both deep stall and perch can be initiated without thrust \cite{taniguchi_AnalysisDeepstall_2008,cunis_DynamicStability_2020,waldock_LearningPerform_2018}, hence, we choose to minimize the thrust used in the landing phase. The landing procedure also needs to avoid an abrupt change in control inputs, so the cost function also encompasses an input deviation term $\Delta \mathbf{u}$. 
\textcolor{black}{
	The pitch constraint in the landing phase, $\theta_r(t)$, for $t_0 \leq t \leq t_f +T_p$, is restricted to be larger than $0^\circ$ to make sure that the UAV will not do a nose-down landing like a normal landing procedure. Especially, the pitch $\theta_r(t)$, for $t_f \leq t \leq t_f +T_p$, is constrained to be $\leq 90^\circ$ to make the perch more natural \cite{venkateswararao_OptimizationStability_2014}.}
Thus, by solving the problem \eqref{eqn_OC_problem_formulation}--\eqref{eqn_OC_problem_formulation_constraints}, we arrive at an optimal trajectory, denoted by the pair of state and input $(\mathbf{x}_r,\mathbf{u}_r)$. This is the reference trajectory used for the tracking mechanism.
%%%%%%%%%%%%%%%%%%%%%%%%%%%%%%%%%%%%%%%%%%%%%%%%%%%%%%%%%%%%%%
\subsection{Landing trajectory tracking}
\label{subsec_Landing_Traject_Track}
Now, we proceed to land the UAV by using NMPC to track the trajectory. The reference optimal trajectory is created without disturbances, notably gusts, from the surrounding environment. The NMPC controller is chosen for the tracking task due to its well-known capability to overcome disturbances, satisfy constraints \cite{grune_NonlinearModel_2017}, \textcolor{black}{and stabilize nonholonomic systems \cite{fontes_GeneralFramework_2001}}.

	Let $\mathbf{x}_e(t) = \mathbf{x}(t) - \mathbf{x}_r(t)$ and $\mathbf{u}_e(t)=\mathbf{u}(t)-\mathbf{u}_r(t)$. We define the error dynamics $\dot{\mathbf{x}}_e = \dot{\mathbf{x}} - \dot{\mathbf{x}}_r = f(\mathbf{x}_e+\mathbf{x}_r,\mathbf{u}_e+\mathbf{u}_r) - f(\mathbf{x}_r,\mathbf{u}_r) \triangleq g(\mathbf{x}_e,\mathbf{u}_e)$, or, in the shortened form
\begin{equation}
\dot{\mathbf{x}}_e = g(\mathbf{x}_e,\mathbf{u}_e).
\end{equation}

The NMPC scheme is stated as follows: at time $t$ $(t_0 \leq t \leq t_f)$, with the error $\mathbf{x}_{e_t}$ (assumed fully measurable), we solve an OCP over the prediction horizon $T_p$
\begin{equation}
	\label{eqn_NMPC_problem_formulation}
	\underset{{\bar{\mathbf{u}}_e(\cdot)}}{\text{min}}   \ J_2\left(\bar{\mathbf{x}}_e, \bar{\mathbf{u}}_e\right)\hspace{-.05cm}
	= \hspace{-.05cm} \int_{t}^{t+T_p} \hspace{-20pt} \ell_2((\bar{\mathbf{x}}_e(\tau), \bar{\mathbf{u}}_e(\tau)) d \tau+ F(\bar{\mathbf{x}}_e(t+T_p)),
\end{equation}
\vspace{-15pt}
\begin{subequations}
	\label{eqn_NMPC_problem_formulation_constraints}
	\begin{align}       
\text{subject to: }		 & \dot{\bar{\mathbf{x}}}_e(\tau) = g(\bar{\mathbf{x}}_e(\tau),\bar{\mathbf{u}}_e(\tau)), \\
		\label{eqn_NMPC_problem_formulation_state_input_constraints}
		&\bar{\mathbf{u}}_e(\tau) + \mathbf{u}_r(\tau) \in \mathcal{U}, \ \tau \in [t,t+T_p] \\
		& \bar{\mathbf{x}}_e(\tau) + \mathbf{x}_r(\tau) \in \mathcal{X}_{\text{MPC}}, \ \tau \in [t,t+T_p]  \\
		\label{eqn_NMPC_problem_formulation_initial_states}
		& \bar{\mathbf{x}}_e(t) = \mathbf{x}_{e_t}, \\
		& \bar{\mathbf{x}}_e(t+T_P) + \mathbf{x}_r(t+T_P) \in \mathcal{X}_\text{MPC}.
	\end{align}
\end{subequations}
Here, the stage cost $\ell_2((\bar{\mathbf{x}}_e(\tau), \bar{\mathbf{u}}_e(\tau))$ and the terminal cost $F(\bar{\mathbf{x}}_e(t+T_p))$ are defined, respectively, by
\begin{equation}
	\label{eqn_NMPC_problem_formulation_costs}
	\| \bar{\mathbf{x}}_e(\tau ) \|_{Q_x}^2\hspace{-5pt} + \| \bar{\mathbf{u}}_e(\tau )\|_{Q_u}^2  ,   \\
	\mbox{and }\|\bar{\mathbf{x}}_e(t+T_p) \|_{Q_{x_f}}^2,  
\end{equation}
the vectors $\bar{\mathbf{x}}_e$, $\bar{\mathbf{u}}_e$ are the predicted state and input errors, $\mathbf{x}_r$, $\mathbf{u}_r$ are the optimal solutions of \eqref{eqn_OC_problem_formulation}--\eqref{eqn_OC_problem_formulation_constraints}, $\{Q_x, {Q_{x_f}} \} \subset \mathbb{S}^6_{+}$, $Q_u \in \mathbb{S}^2_{++}$ are weighting matrices, $\mathcal{U}$ is the same input constraint \eqref{eqn_OC_problem_formulation_state_input_constraints} as in the generation task, and $\mathcal{X}_{\text{MPC}}$ is  the state constraint set as in \eqref{eqn_set_X_MPC}, which serves as the stage constraint set and as the terminal constraint set. After obtaining the optimal solution $\mathbf{u}(\cdot)$, its values in the period $[t,t+\delta t]$, with $\delta t < T_p$, are applied to the system. At $t+\delta t $, the state is sampled, time shifts $ t \leftarrow t+\delta t$, and the procedure runs recurrently from $t_0$ to $t_f$.
\textcolor{black}{	In the tracking task, the pitch $\theta$ and the pitch rate $\omega_y$ are relaxed to $\theta \in [\underaccent{\bar}{\theta}_\text{MPC},\bar{\theta}_\text{MPC}], 
	|\omega_y| \leq \bar{\omega}_{y_\text{MPC}} $ to provide some flexibility for the FW UAV in windy conditions.}

%%%%%%%%%%%%%%%%%%%%%%%%%%%%%%%%%%%%%%%%%%%%%%%%%%%%%%%%%%%%%%
%%%%%%%%%%%%%%%%%%%%%%%%%%%%%%%%%%%%%%%%%%%%%%%%%%%%%%%%%%%%%%
%%%%%%%%%%%%%%%%%%%%%%%%%%%%%%%%%%%%%%%%%%%%%%%%%%%%%%%%%%%%%%
\section{SIMULATION}
\label{sec_Simulation}
%%%%%%%%%%%%%%%%%%%%%%%%%%%%%%%%%%%%%%%%%%%%%%%%%%%%%%%%%%%%%%
%%%%%%%%%%%%%%%%%%%%%%%%%%%%%%%%%%%%%%%%%%%%%%%%%%%%%%%%%%%%%%
\subsection{Aerosonde UAV}
\label{subsec_Aerosonde_UAV}
The FW UAV model implemented in this paper is the Aerosonde FW UAV, whose physical and aerodynamic parameters are provided in \cite[Appendix E.2]{beard_SmallUnmanned_2012}. 
\begin{figure}[ht!]
	\begin{center}
		\includegraphics[width=0.7\columnwidth]{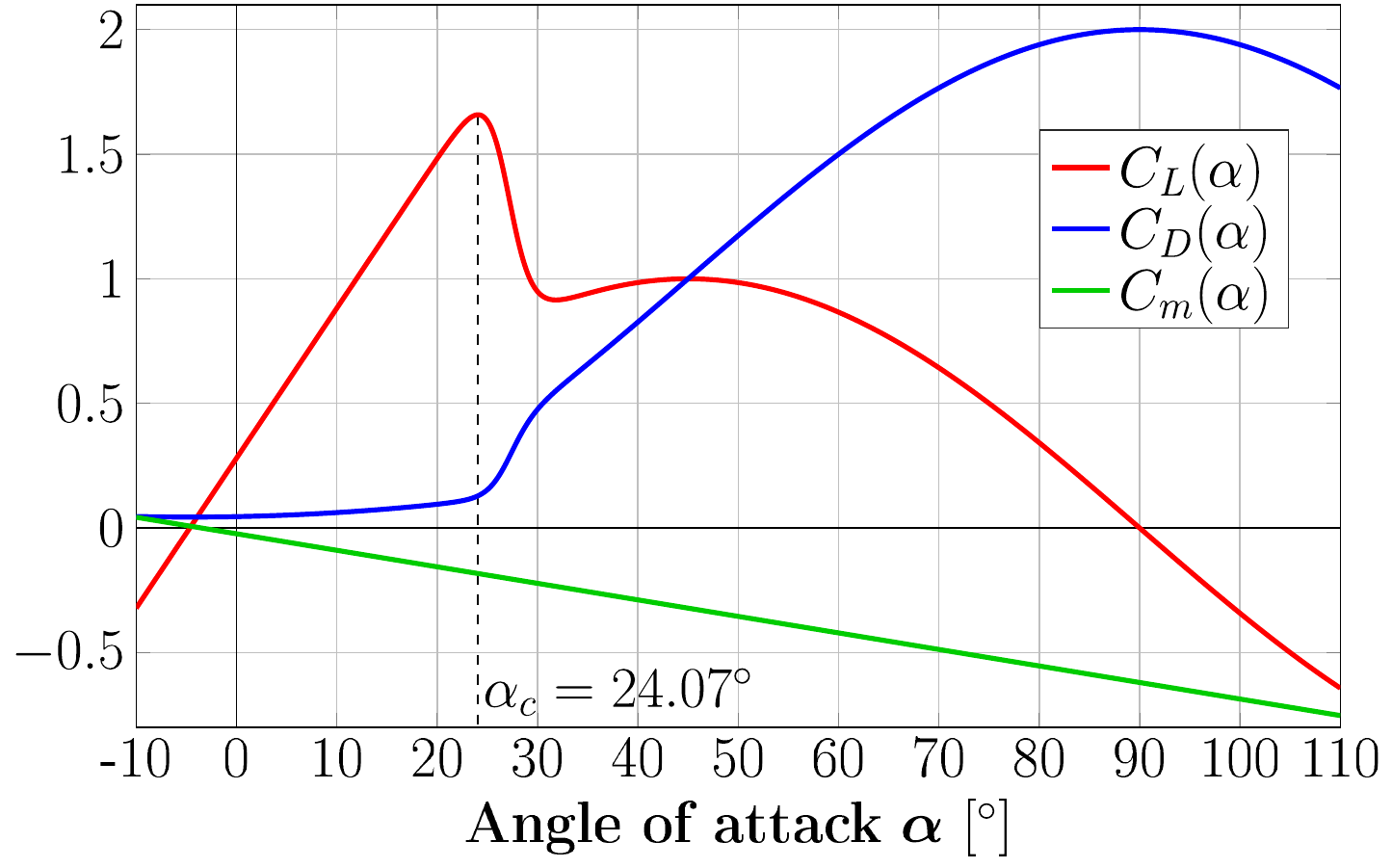}
		\caption{Aerodynamic coefficients.}
		\label{fig_NACA_4412_AR_12_coefficients}
	\end{center}
\end{figure}
The aerodynamic coefficients in \eqref{eqn_C_L_C_D_C_m} are plotted in Fig. \ref{fig_NACA_4412_AR_12_coefficients} for $\alpha \in [-10^\circ,110^\circ]$, whose critical AoA value is $\alpha_c = 24.07^\circ$.
%%%%%%%%%%%%%%%%%%%%%%%%%%%%%%%%%%%%%%%%%%%%%%%%%%%%%%%%%%%%%%
%%%%%%%%%%%%%%%%%%%%%%%%%%%%%%%%%%%%%%%%%%%%%%%%%%%%%%%%%%%%%%
\subsection{Dryden wind turbulence model}
\label{subsec_Dryden}
To make the simulation as realistic as possible, the wind components $u_w$, $w_w$ in \eqref{eqn_state_3_state_4_wind} contain not only the known steady winds, denoted with the subscript $\square_s$, but also random gusts, \textcolor{black}{denoted by $\square_g$}. 
\begin{equation}
	\label{eqn_wind_gust}
	\begin{bmatrix}
		u_w \\
		w_w
	\end{bmatrix}
	=
	\mathcal{R}^{\mathcal{B}}_{\mathcal{I}}(\theta)
	\begin{bmatrix}
		{u}_s^\mathcal{I} \\
		{w}_s^\mathcal{I}
	\end{bmatrix}
	+
	\begin{bmatrix}
		{u}_g^\mathcal{B} \\
		{w}_g^\mathcal{B} 
	\end{bmatrix}
\end{equation}
The gusts follow the Dryden wind turbulence model \cite{beard_SmallUnmanned_2012}, which are generated by passing the white noise signals, denoted by $\Omega_u^\mathcal{B}$ and $\Omega_w^\mathcal{B}$, through the transfer functions
%\begin{subequations}
%\begin{align}
% H_{u}(s) &=\sigma_{u} \sqrt{\frac{2 V_{a_0} }{L_{u}}} \frac{1}{s+\frac{V_{a_0}}{L_{u}}}, \\
%H_{w}(s) &=\sigma_{w} \sqrt{\frac{3 V_{a_0} }{L_{w}}} \frac{\left(s+\frac{V_{a_0}}{\sqrt{3} L_{w}}\right)}{\left(s+\frac{V_{a_0}}{L_{w}}\right)^{2}},
%\end{align}
%\end{subequations}
\begin{subequations}
	\begin{gather}
		H_{u}(s) =\sigma_{u} \sqrt{ \tfrac{2 V_{a_D}}{\pi L_{u}}}  \tfrac{ 1}{s+ V_{a_0} / L_{u}}, \\
		H_{w}(s) =\sigma_{w} \sqrt{ \tfrac{3 V_{a_D}}{\pi L_{w}}} \tfrac{1}{\left(s+V_{a_0}/L_{w}\right)^{2}},
	\end{gather}
\end{subequations}
in which $\sigma_u,\sigma_w \ [m/s]$ are the turbulence intensities along the body frame $x$ and $z$ axes, $L_u,L_w \ [m]$ are spatial wavelengths, $V_{a_D} \ [m/s]$ is the airspeed of the aircraft assumed to be constant. For the simulation, we choose the ``low altitude, light turbulence'' scenario with $\sigma_u=1.06$, $\sigma_w=0.7$, $L_u=200$, $L_w=50$, $V_{a_D}=25$. We choose the noise $\Omega_u^\mathcal{B}$ $\sim \mathcal{N}(0,0.5)$, $\Omega_w^\mathcal{B} $ $\sim \mathcal{N}(0,0.5)$, and bound the gusts $u_g^\mathcal{B}, w_g^\mathcal{B} \in [-0.2,0.2] (m/s)$. The NMPC is solved by taking into account the steady wind components, but without information about the gusts. 
%%%%%%%%%%%%%%%%%%%%%%%%%%%%%%%%%%%%%%%%%%%%%%%%%%%%%%%%%%%%%%
%%%%%%%%%%%%%%%%%%%%%%%%%%%%%%%%%%%%%%%%%%%%%%%%%%%%%%%%%%%%%%
\subsection{Simulation parameters}
\label{subsec_Simulation_Params}
The specific input constraint set $\mathcal{U}$ and the state constraint sets $\mathcal{X}_{g_1}$, $\mathcal{X}_{g_2}$, $\mathcal{X}_{f}$, and $\mathcal{X}_\text{MPC}$ in \eqref{eqn_set_constraints} are gathered in Table \ref{table_problem_formulation}. 
$c_1$ in \eqref{eqn_req_track_Va} is chosen to be $0.3$ \cite{feroskhan_DynamicsSideslip_2016}, we choose $c_2$ in \eqref{eqn_req_gen_Va} to be $0.25$, $\delta x$ in \eqref{eqn_req_gen_xf} to be $10^{-4} \ m$.
The size of the recovery net is chosen to be $5m$ width $\times$ $3m$ height. Since the Aerosonde UAV has the length of $1.7m$ \cite{burston_ReverseEngineering_2014}, $\bar{z}_\text{net} = -1.7m$ and $\underaccent{\bar}{z}_\text{net}=-4.7m$.
The recovery is deemed successful when the center of mass of the UAV lands inside the surface of the net, hence, in the generation task, we choose $\Delta z=0.5m$ (in \eqref{eqn_req_gen_zf}), which makes the reference $z_r(t)$ for $t_f \leq t \leq t_f +T_p$ to be limited to $-4.2m$ to $-2.2m$.
Before arriving at the net, to make sure that the UAV does not slam into the ground, in the \emph{generation task}, for $t_0\leq t < t_f$, the altitude is constrained to be $\geq 2m$ ($z_r(t) \leq -2m$, $\delta_z$ in \eqref{eqn_set_X_g_1} is $0.3m$). The pitch rate is restricted in $-1.46 \leq \omega_{y_r}(t)  \leq 1.46  $. 
In the NMPC \emph{tracking task}, the constraints on pitch and pitch rate are relaxed to $-50^\circ \leq \theta(t) \leq 180^\circ$ and $-\pi/2 \leq \omega_y(t) \leq \pi/2 $.

From \eqref{eqn_airspeed_and_AoA} and \eqref{eqn_wind_gust}, when there is no gust, $u$ and $w$ are calculated as functions of the airspeed $V_a$, the AoA $\alpha$, and the steady winds $u_s^\mathcal{I}$, $w_s^\mathcal{I}$ as follows
		\begin{equation}
				\label{eqn_vx_vz_aoa_va}
		\begin{bmatrix}
			u \\
			w
		\end{bmatrix}
		= V_a
			\begin{bmatrix}
			\cos{\alpha} \\
			\sin{\alpha}
		\end{bmatrix}
	+
		\mathcal{R}^{\mathcal{B}}_{\mathcal{I}}(\theta)
		\begin{bmatrix}
			u_s^\mathcal{I} \\
			w_s^\mathcal{I}
		\end{bmatrix}.
	\end{equation}
%\begin{equation}
%	\label{eqn_vx_vz_aoa_va}
%	u = V_a \cos{\alpha} + u_w, \ 
%	w = V_a \sin{\alpha} + w_w.
% \end{equation}
By using the relations \eqref{eqn_vx_vz_aoa_va}, we can choose the values of $u$ and $w$ in \eqref{eqn_uav_2d_dynamics} through choosing the values of the AoA $\alpha$, the airspeed $V_a$, and the steady wind component in $\mathcal{I}$. Thus, the state vector $\mathbf{x}$ can be expressed as a function of another vector $\bm{\xi}$, $\mathbf{x}=\kappa(\bm{\xi})$, where $\bm{\xi}  \triangleq \left[ x,z,V_a,\alpha,\theta,\omega_y, u_s^\mathcal{I},w_s^\mathcal{I} \right]^\intercal$, to facilitate the choice of the boundary constraints. 

For the boundary constraints, in the \emph{generation and tracking tasks}, at $t_0$ we impose the following specific constraint which obtained by solving \eqref{eqn_level_trim_flight} for $V_{a_r}=25m/s$, $\bm{x}_r(t_0)=\kappa(\bm{\xi}_{0_r})$, $ \bar{\mathbf{x}}_e(t_0) = \bm{0} $, with 
\begin{equation}
	\label{eqn_initial_boundary}
 \bm{\xi}_{0_r} = [-280 ,-200,25,\alpha_0,\theta_0,0,u_s^\mathcal{I},0]^\intercal. 
\end{equation}

\begin{table*}[b]
	\caption{State and input constraints}
	\label{table_problem_formulation}
	\begin{center}
		\begin{tabular}{c | c | c | c |c | c | c }
			\multirow{2}{*}{\textbf{Parameters}} & \multicolumn{4}{c|}{\textbf{Generation task}} & \multicolumn{2}{c}{\textbf{NMPC trajectory tracking}} \\
			\cline{2-7}
			& $\bm{t_0}$          
			&  $\bm{t_0 < t < t_f}$   
			& $\bm{t_f}$ 
			& $\bm{t_f < t \leq (t_f + T_p)}$ 
			& $\bm{t_0}$ 
			& $\bm{t_0 < t \leq t_f}$\\ 
			\hline
			\makecell{State \\ constraint set} & \multicolumn{2}{c|}{$\mathcal{X}_{g_1}$} & $\mathcal{X}_{f} \subset \mathcal{X}_{g_2}$ & $\mathcal{X}_{g_2}$& \multicolumn{2}{c}{$\mathcal{X}_\text{MPC}$} \\
			\hline
			$x[m]$ & $x_0$ & $x(t) \leq 0$ & $ 0 \leq x(t_f) \leq 10^{-4}$ &$  x(t) \geq  0$ &  $x_0$ &\multirow{5}{*}{\resizebox{0.27\textwidth}{!}{\makecell{$
					\mathcal{C}_{\text{pos}} = \left\{ x,z \middle|
					\begin{aligned}
						-z -4.7 +0.9812 x \leq & 0 \ \text{if} \ x<0\\
						-z -  4.7 \leq & 0 \ \text{if} \ x \geq 0\\
						z + 1.7 \leq  & 0 
					\end{aligned} \right\}
					$ \\
					$
					\mathcal{C}_{V_a} = \left\{ x,V_a \middle|
					\begin{aligned}
						V_a - 7.5 + 0.55 x \leq & 0 \ \text{if} \ x < 0 \\
						V_a - 7.5 \leq &  0 \ \text{if} \ x \geq 0
					\end{aligned} \right\}
					$ \\		}	
				 }} \\
			\cline{1-6}
			$z[m]$ & $z_0$  & $ z(t)\leq -2 $ & \multicolumn{2}{c|}{$ -4.2 \leq z(t) \leq -2.2 $}  & $z_0$ \\
			\cline{1-6}
			$u[m/s]$ & $u_0$  & $-40 \leq u(t)\leq 40 $ &  \multicolumn{2}{c|}{--}  & $u_0$ \\
			\cline{1-6}
			$w[m/s]$ & $w_0$ & $-40 \leq w(t)\leq 40 $ &\multicolumn{2}{c|}{--}  & $w_0$\\
			\cline{1-6}
			$V_a[m/s]$ & $V_{a_0}$ & -- & \multicolumn{2}{c|}{$V_a(t) \leq 7$}   & $V_{a_0}$ \\
			\cline{1-6}
			$\alpha^\circ$ & $\alpha_0$ &$-10^\circ \leq \alpha(t) \leq 110^\circ$ & \multicolumn{2}{c|}{$90^\circ \leq \alpha(t) \leq 110^\circ $}  & $\alpha_0$ & 	$-10^\circ \leq \alpha(t) \leq 110^\circ$ \\
			\hline
			$\theta^\circ$  & $\theta_0$ & $0^\circ \leq \theta(t) \leq 150^\circ$ &  \multicolumn{2}{c|}{$0^\circ \leq \theta(t) \leq 90^\circ$ }   &$\theta_0$ & $-50^\circ \leq \theta(t) \leq 180^\circ$ \\
			\hline
			$\omega_y [rad/s]$  & $0 $ & \multicolumn{3}{c|}{$-1.46 \leq \omega_y(t)   \leq 1.46 $}& $0 $ &  $-\pi/2 \leq \omega_y(t)   \leq \pi/2 $\\
			\hline
			$\mathbf{u}(t) \in \mathcal{U}$ & $\bm{u}_0 $& \multicolumn{3}{c|}{$\left[ -40^\circ , 0 \right]^\intercal \leq \left[ \delta_{e} , \delta_{t} \right]^\intercal \leq \left[40^\circ , 1 \right]^\intercal $ } & -- &$\left[ -40^\circ , 0 \right]^\intercal \leq \left[ \delta_{e} , \delta_{t} \right]^\intercal \leq \left[40^\circ , 1 \right]^\intercal $ \\
		\end{tabular}   
	\end{center}
\end{table*}

We choose $t_f$ \textcolor{black}{as in \eqref{eqn_OC_problem_formulation}--\eqref{eqn_OC_problem_formulation_constraints}} to be $24s$, the prediction horizon $T_p$ \textcolor{black}{as in \eqref{eqn_OC_problem_formulation}--\eqref{eqn_OC_problem_formulation_constraints},  \eqref{eqn_NMPC_problem_formulation}--\eqref{eqn_NMPC_problem_formulation_constraints}} to be $0.5s$ \textcolor{black}{and the sampling time $\delta t$ in \eqref{eqn_req_track_xf} and in Section \ref{subsec_Landing_Traject_Track} is chosen as $0.1s$}.
The OCPs \eqref{eqn_OC_problem_formulation}--\eqref{eqn_OC_problem_formulation_constraints} and \eqref{eqn_NMPC_problem_formulation}--\eqref{eqn_NMPC_problem_formulation_constraints} are transcribed into Nonlinear Programming Problems (NLPs) by using the ``Direct Multiple Shooting'' method \cite{bock_MultipleShooting_1984}. 
For the OCP \eqref{eqn_OC_problem_formulation}--\eqref{eqn_OC_problem_formulation_constraints}, in $[t_0, t_f+T_p]$, states and control inputs are discretized into $N=245$ arcs, equivalent to $t_0,t_1,\dots,t_N=t_f+T_p$. 
In each arc, $t\in[t_i,t_{i+1}]$ \textcolor{black}{$(i \in \{0,1,\dots,N-1\})$}, the inputs and states are parametrized as decision variables, where the input is \textcolor{black}{kept constant} and the system dynamics \eqref{eqn_uav_2d_dynamics} is solved with an arbitrary initial value. The solution of the ODE \eqref{eqn_uav_2d_dynamics} at time $t_{i+1}$ is obtained with the Runge--Kutta $4^{th}$--order algorithm, with discretization step $\delta t=0.1$, being each state constrained to be equal to the initial value of the next arc. 
The same procedure is applied for the OCP \eqref{eqn_NMPC_problem_formulation}--\eqref{eqn_NMPC_problem_formulation_constraints}. In each NMPC iteration $[t, t+T_p]$, states and control inputs are discretized into $N'=5$ arcs, equivalent to $t,t+\delta,\dots,t_{N'}=t+5\delta =t+T_P$. However, only the input in the first arc $[t, t+\delta]$ is applied to the system.
The NLPs are then solved by using the interior point method in the IPOPT solver \cite{wachter_ImplementationInteriorpoint_2006} within the CasADi toolbox \cite{andersson_CasADiSoftware_2019}.
For the OCP \eqref{eqn_OC_problem_formulation}--\eqref{eqn_OC_problem_formulation_constraints}, the initial guesses for the states are linearly interpolated between 2 points 
$\bm{\xi}(t_0) = [-280, \allowbreak-200,\allowbreak 25,\allowbreak \alpha_0,\allowbreak \theta_0,\allowbreak0,\allowbreak u_s^\mathcal{I},\allowbreak0]^\intercal 
\text{and} \ 
\bm{\xi}(t_f+T_p)= [0,\allowbreak-3.2,\allowbreak7,\allowbreak 110^\circ,\allowbreak 90^\circ,\allowbreak0,\allowbreak u_s^\mathcal{I},\allowbreak0]^\intercal$, 
while the initial guesses for the inputs are $[0,0]^\intercal$. 
For the OCP \eqref{eqn_NMPC_problem_formulation}--\eqref{eqn_NMPC_problem_formulation_constraints}, the initial guesses are $(\mathbf{x}_r,\mathbf{u}_r)$. 

In \eqref{eqn_OC_generation_cost}, the weighting terms are chosen as follows
${P_{\delta_t}}= 1000$, ${P_{\Delta U}} = 4000  \mathbb{I}_2$. In 
\eqref{eqn_NMPC_problem_formulation_costs}, the weighting terms are ${Q_x}=\text{diag}\{200,200,10,10,1,1\}$, ${Q_{x_f}} = 10{Q_x}$, ${Q_u}=20\mathbb{I}_2$. Simulations are run on a lab computer with an AMD Ryzen $5$ $2600$ $6$-core processor, $3.4$$GHz$, $12$ CPUs, $16GB$ RAM, Python $3.8.8$.

%%%%%%%%%%%%%%%%%%%%%%%%%%%%%%%%%%%%%%%%%%%%%%%
\subsection{\textcolor{black}{Results and }Analysis}
\label{subsec_result_analysis}
\begin{figure}[htb]
	\begin{center}		 	
		\includegraphics[width=0.9\columnwidth]{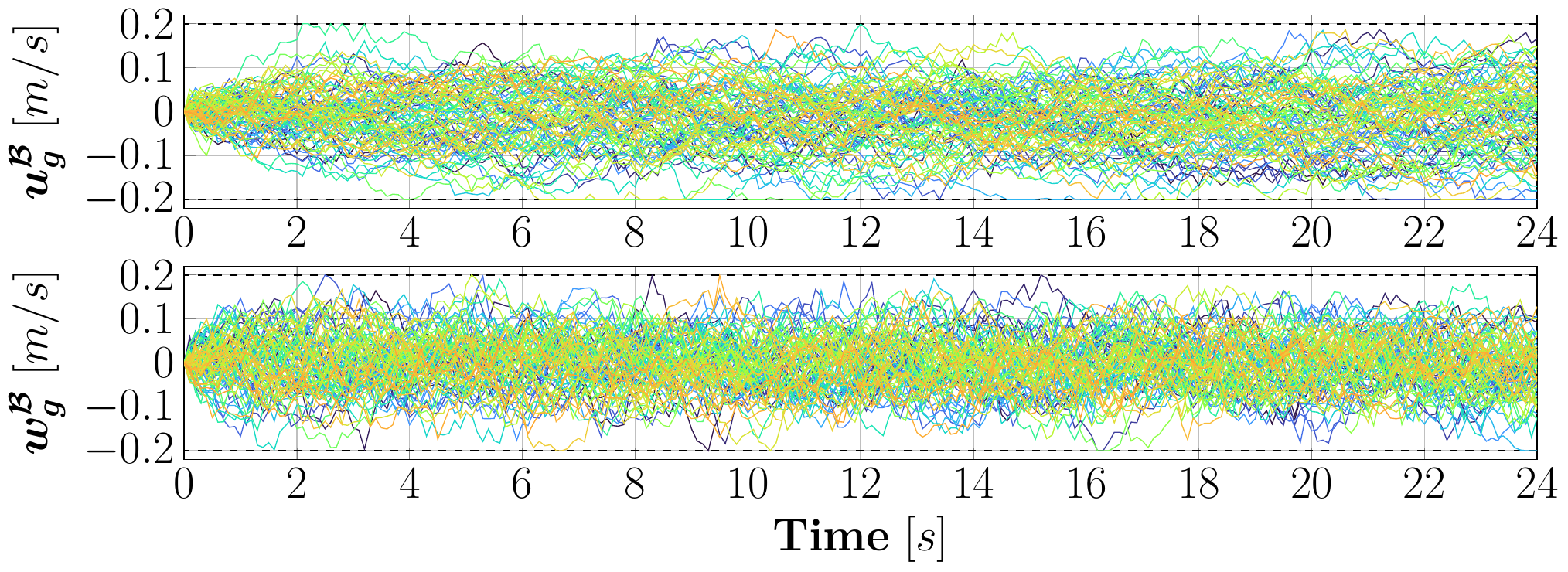}
		%\vspace{-2pt}
		\caption{Gust in the body frame for the $90$ tracking cases.}
		\label{fig_gust_body_frame_satisfied}
	\end{center}
\end{figure}

We solve the OCP \eqref{eqn_OC_problem_formulation}--\eqref{eqn_OC_problem_formulation_constraints} for constant wind $u_s^\mathcal{I}$ from $2m/s$ (tailwind) to $-6m/s$ (headwind), with $1m/s$ step, and vary the time in \eqref{eqn_OC_problem_formulation} to be $t \in [t_0,t_f + T_p + (-10)u_s^\mathcal{I}]$, i.e., the stronger the headwind, the longer the simulation time for the generation task. Only the headwind from $-6m/s$ to $-1m/s$ \textcolor{black}{and the nominal scenario $u_s^\mathcal{I}=0m/s$} give feasible solutions, and their trajectories are plotted in Fig. \ref{fig_pos_with_Va_NMPC}, while the AoA and the airspeed are in Fig. \ref{fig_alpha_airspeed_multiple_wind} in Appendix \ref{sec_appendix_wind}. 
Then, we choose to track the trajectory with $u_s^\mathcal{I}=0m/s$ for $100$ times \textcolor{black}{with the bounded Dryden gust in Subsection \ref{subsec_Dryden}}, only $92$ successful tracking results, and $90$ scenarios satisfy the control requirements \eqref{eqn_requirements}, $2$ scenarios violated the airspeed constraint \textcolor{black}{$V_a(t_f) \leq 7.5m/s$}. \textcolor{black}{The reason is because the unknown gusts to the NMPC tracking controller ``push'' the airspeed at $t_f$ out of the LoS corridor $\mathcal{C}_{V_a}$}. In these $90$ cases, the mean and standard deviation of the final $x$, $z$, and $V_a$ are in Table \ref{table_trajectory_tracking_results}.
\begin{table}[ht]
	\caption{Trajectory tracking results}
	\label{table_trajectory_tracking_results}
	\begin{center}
		\begin{tabular}{cccc}
		           & $x(t_f)$& $z(t_f)$ & $V_a(t_f)$  \\
			\hline
			Mean $[m]$   & $0.0399$               & $-4.2292$            & $7.0957$      \\
			\hline
			Standard deviation $[m]$  & $0.0476$                 & $0.0250$              & $0.0967$    \\
		\end{tabular}    
	\end{center}
\end{table}
The bounded Dryden gusts in Section \ref{subsec_Dryden} for the $90$ successful tracking scenarios are shown in Fig. \ref{fig_gust_body_frame_satisfied}.
% wind + gust 

\begin{figure}[htb]
	\begin{center}		 	
		\includegraphics[width=0.9\columnwidth]{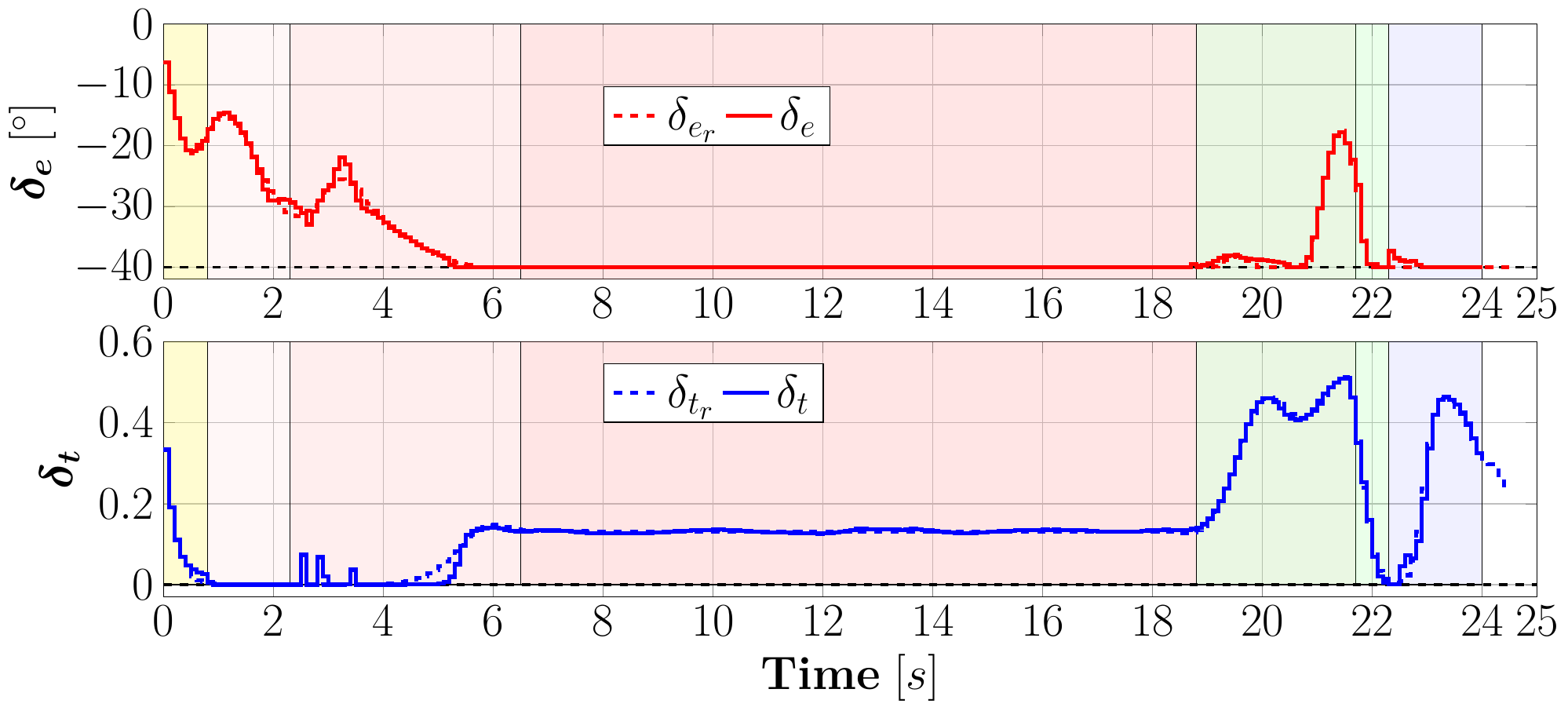}	
		%\vspace{-2pt}
		\caption{Control inputs from NMPC tracking.}
		\label{fig_input_NMPC}
	\end{center}
\end{figure}

\begin{figure}[htb]
	\begin{center}		 	
		\includegraphics[width=0.9\columnwidth]{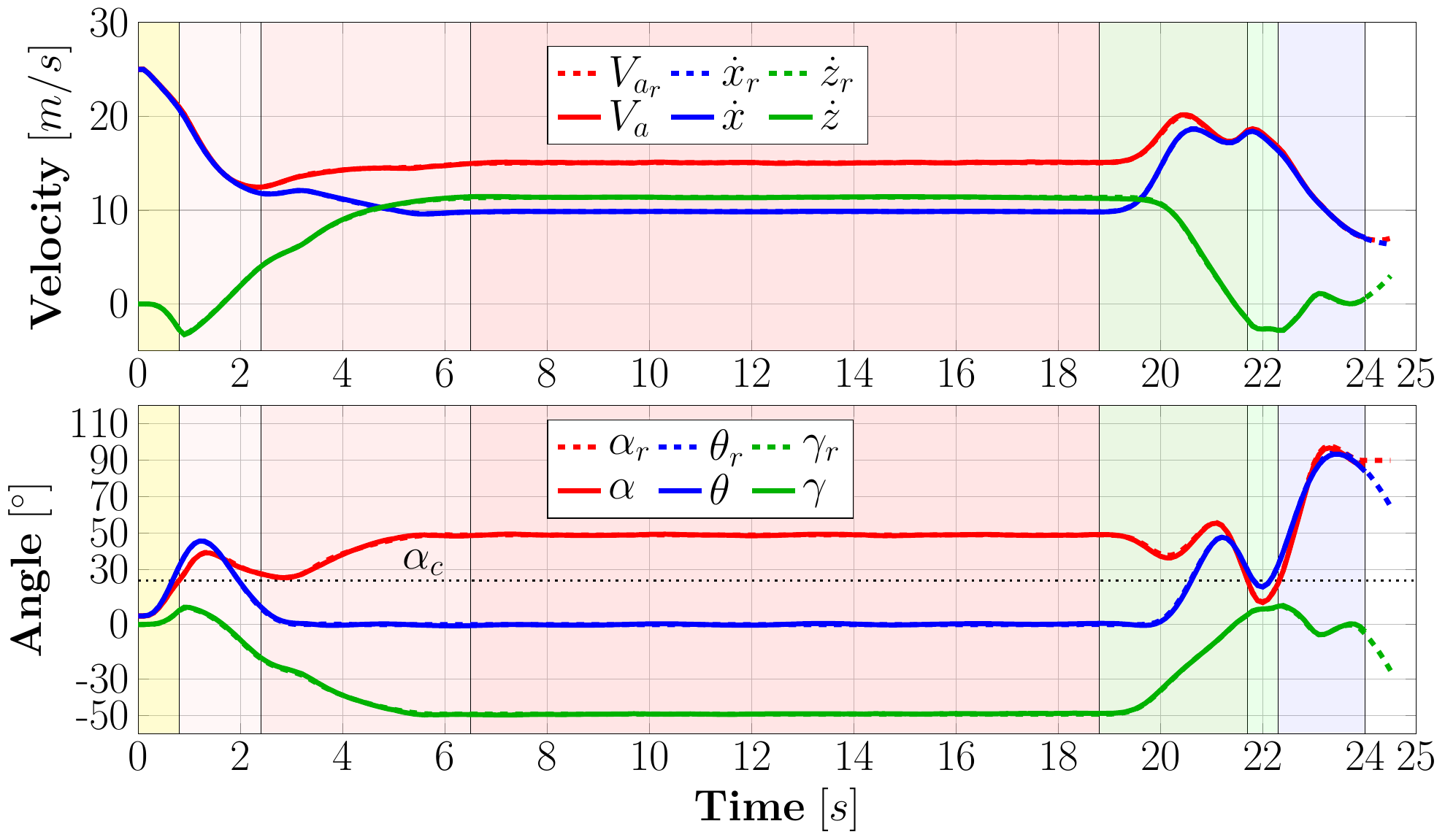}	
		%\vspace{-2pt}
		\caption{Velocities and angles from NMPC tracking.}
		\label{fig_velocity_angle_NMPC}
	\end{center}
\end{figure}

\begin{figure}[htb]
	\begin{center}		 	
		\includegraphics[width=0.9\columnwidth]{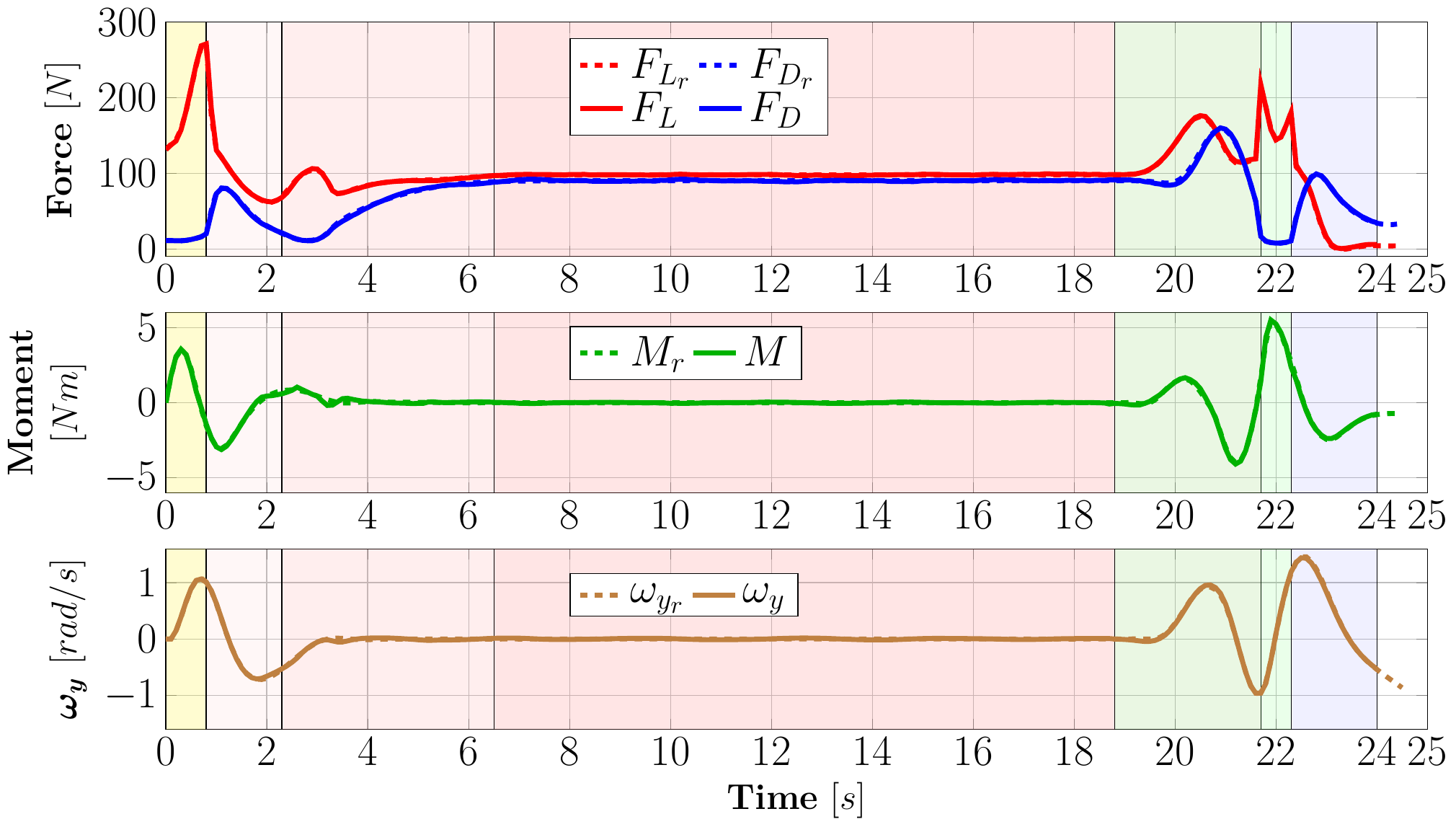}
		%\vspace{-2pt}
		\caption{Aerodynamics and pitch rate from NMPC tracking.}
		\label{fig_aerodynamics_NMPC}
	\end{center}
\end{figure}

\textcolor{black}{We plot one successful NMPC tracking result in Figs. \ref{fig_pos_with_Va_NMPC}--\ref{fig_aerodynamics_NMPC}, where the reference values obtained from \eqref{eqn_OC_problem_formulation}--\eqref{eqn_OC_problem_formulation_constraints} are plotted in dashed lines, and the trajectory tracking results of \eqref{eqn_NMPC_problem_formulation}--\eqref{eqn_NMPC_problem_formulation_constraints} are plotted in solid lines.}
Fig. \ref{fig_pos_with_Va_NMPC} shows the landing trajectory tracking results. The green images of the UAV are plotted every $0.3s$ to show the position in the $xz$ plane, orientation ($\theta$), and elevator deflection angle ($\delta_e$). The length of the UAV is enlarged to $10m$ solely for visualization. 
% position 
\begin{figure}[htb]
	\begin{center}		 	
		\includegraphics[width=0.8\linewidth]{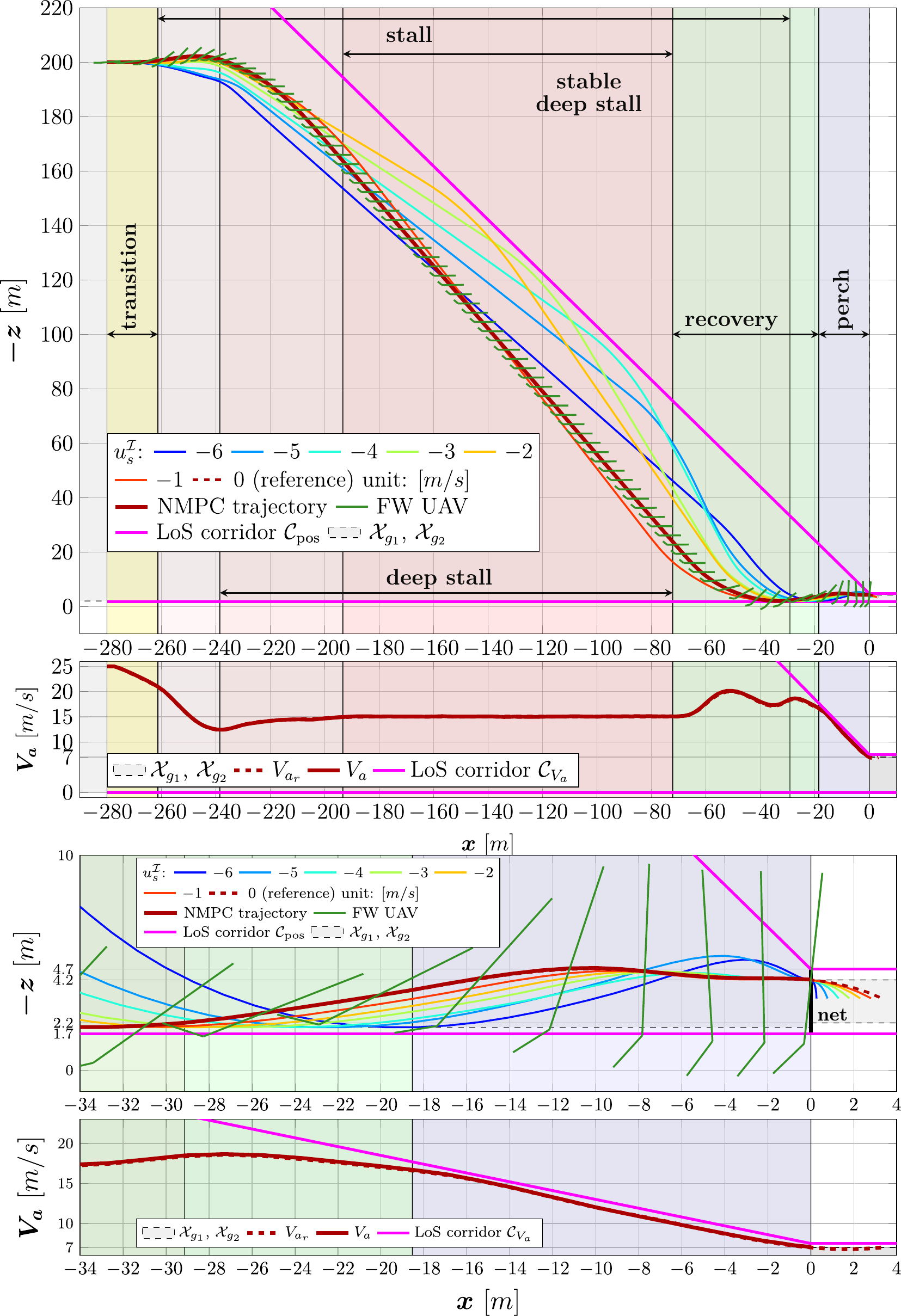}	
		%\vspace{-8pt}
		\caption{Trajectory tracking results.}
		\label{fig_pos_with_Va_NMPC}
	\end{center}
\end{figure}
The landing trajectory is composed of four phases: transition from cruise to stall, stall to deep stall, deep stall recovery, and perch to land on the net. These phases are presented in light yellow, red, green, and blue background colors in Figures \ref{fig_pos_with_Va_NMPC} to \ref{fig_aerodynamics_NMPC}.

The UAV prepares its attitude for stall by quickly pitching up (Fig. \ref{fig_velocity_angle_NMPC}), increases its AoA, which makes it gain some altitude, and then eventually fall into the stall state. This maneuver matches the transition in \cite{taniguchi_AnalysisDeepstall_2008,cunis_DynamicStability_2020,mathisen_PrecisionDeepStall_2021}. 
As soon as the UAV enters stall at $0.8s$ (the AoA surpasses its critical value $\alpha_c$), there is a sudden drop in lift, an increase in drag (Fig. \ref{fig_aerodynamics_NMPC}), and the airspeed is decreased (Fig. \ref{fig_velocity_angle_NMPC}), while the thrust is retained as $0$ (Fig. \ref{fig_input_NMPC}). After the period of $1.5s$, the UAV reaches deep stall (medium red background), the airspeed and lift start to regain. 
From $6.5s$ to $18.8s$, the airflow over the wings of the aircraft becomes stable, and the aircraft falls into the stable deep stall state, which is emphasized by the dark red background in Figs. \ref{fig_pos_with_Va_NMPC}--\ref{fig_aerodynamics_NMPC}. The airspeed and AoA are steady in this phase, even in the windy condition. 
At the end of deep stall $(18.8s)$, the UAV goes through the recovery phase, where it decreases its AoA to return to the normal operating region. The recovery phase could either be commenced by pitching down or by increasing the thrust \cite{faa_AirplaneFlying_2021,cunis_DynamicStability_2020}. Since we impose $\theta_r(t) \geq 0$ in the landing phase of the generation task, the UAV must increase its thrust at $18.8s$. Even though the airspeed is already high at this moment $(>15m/s)$, this recovery technique feeds more speed for the UAV and it is counter-intuitive that the large kinetic energy at the end of deep stall could be mitigated by transforming into the potential energy when perch. 
To compensate for the increase in thrust, at around $20s$, the pitch is increased to generate more drag, which slows down the UAV. 
The combination of increasing the airspeed and lessening the AoA leads to a surge in the lift (Fig. \ref{fig_aerodynamics_NMPC}) that also prevents the aircraft from slamming into the ground.  
When the UAV escapes stall (at $21.7s$), the lift is increased and the drag is decreased, as opposed to the beginning of stall at $0.8s$.
Finally, in the perching phase, the UAV pitches up, increasing its AoA and drag. However, the thrust is employed and the UAV is almost perpendicular to the ground to hold its altitude.
% input
\vspace{-15pt}
% velocity 
\vspace{-15pt}
% aerodynamics 

%  \ref{table_comparison_landing_results} compares our results with previous work in literature, with $\Delta h = h(t_f)-h(t_l)$, $\Delta x = x(t_f)-x(t_l)$, and $\Delta T = t_f-t_l$, which is in accordance with Table \ref{table_comparison_landing_techniques}. In \cite{mathisen_PrecisionDeepStall_2021}, the UAV has to climb up for deep stall and travels a longer distance, hence, their landing time is significantly higher.
%%%%%%%%%%%%%%%%%%%%%%%%%%%%%%%%%%%%%%%%%%%%%
\begin{comment}
\begin{table*}[t]
	\caption{Comparison of landing results for fixed-wing UAVs}
	\label{table_comparison_landing_results}
	\begin{center}
		\begin{tabular}{ccccccc}
			\textbf{Landing technique}  &$m$ $[kg]$ &$\Delta h$ $[m]$ & $\Delta x$ $[m]$ & $V_a(t_l)$ $[m/s]$ & $V_a(t_f)$ $[m/s]$ & $\Delta T$ $[s]$  \\
			\hline
			\makecell{Deep stall, $2 m/s$ wind  \cite{mathisen_PrecisionDeepStall_2021}}  & $13.5$ & $-200$  & $340$ & $25$ & $9.67$            & $48$      \\
			\hline
			\makecell{Perch  \cite{venkateswararao_ParametricStudy_2015}}  &$0.8$ & $5.3$   & $48$    & $14.6$       & $0.7$              & $4.3$    \\
			\hline
			Our approach, $2m/s$ wind &  $13.5$ & $-100$     & $177$   &  $21.95$      &   $9.18$         & $13.5$   
		\end{tabular}    
	\end{center}
\end{table*}
\end{comment}
%%%%%%%%%%%%%%%%%%%%%%%%%%%%%%%%%%%%%%%%%%%%%%%%%%%%%%%%%%%%%%
\section{CONCLUSIONS}
\label{sec_Conclusions}
This paper successfully demonstrates a new recovery technique for FW UAVs by combining deep stall with perch.
Future work will focus on the following key challenges: (i) increase the versatility of the initialization of optimization problems to guarantee their feasibility; (ii) investigate novel robust NMPC schemes adapted to this particular problem; and (iii) design more efficient computational implementations to reduce the time generating the landing trajectory offline and the on-line tracking controls.

%%%%%%%%%%%%%%%%%%%%%%%%%%%%%%%%%%%%%%%%%%%%%%%%%%%%%%%%%%%%%
%\appendix{Various trajectories based on different initial wind conditions}

%Here, we 

%\appendix{something else}
%%%%%%%%%%%%%%%%%%%%%%%%%%%%%%%%%%%%%%%%%%%%%%%%%%%%%%%%%%%%%

%%%%%%%%%%%%%%%%%%%%%%%%%%%%%%%%%%%%%%%%%%%%%%%%%%%%%%%%%%%%%
% \section{ACKNOWLEDGMENTS}
% \label{sec_Acknowledgments}
\vspace{-4pt}
%%%%%%%%%%%%%%%%%%%%%%%%%%%%%%%%%%%%%%%%%%%%%%%%%%%%%%%%%%%%%
%%%%%%%%%%%%%%%%%%%%%%%%%%%%%%%%%%%%%%%%%%%%%%%%%%%%%%%%%%%%%
\begin{appendices}
	%%%%%%%%%%%%%%%%%%%%%%%%%%%%%%%%%%%%%%%%%%%%%%%%%%%%%%%%%%%%%
	\section{Various initial wind conditions}
	\label{sec_appendix_wind}
\begin{figure}[htb]
		\begin{center}		 	
			\includegraphics[width=0.8\columnwidth]{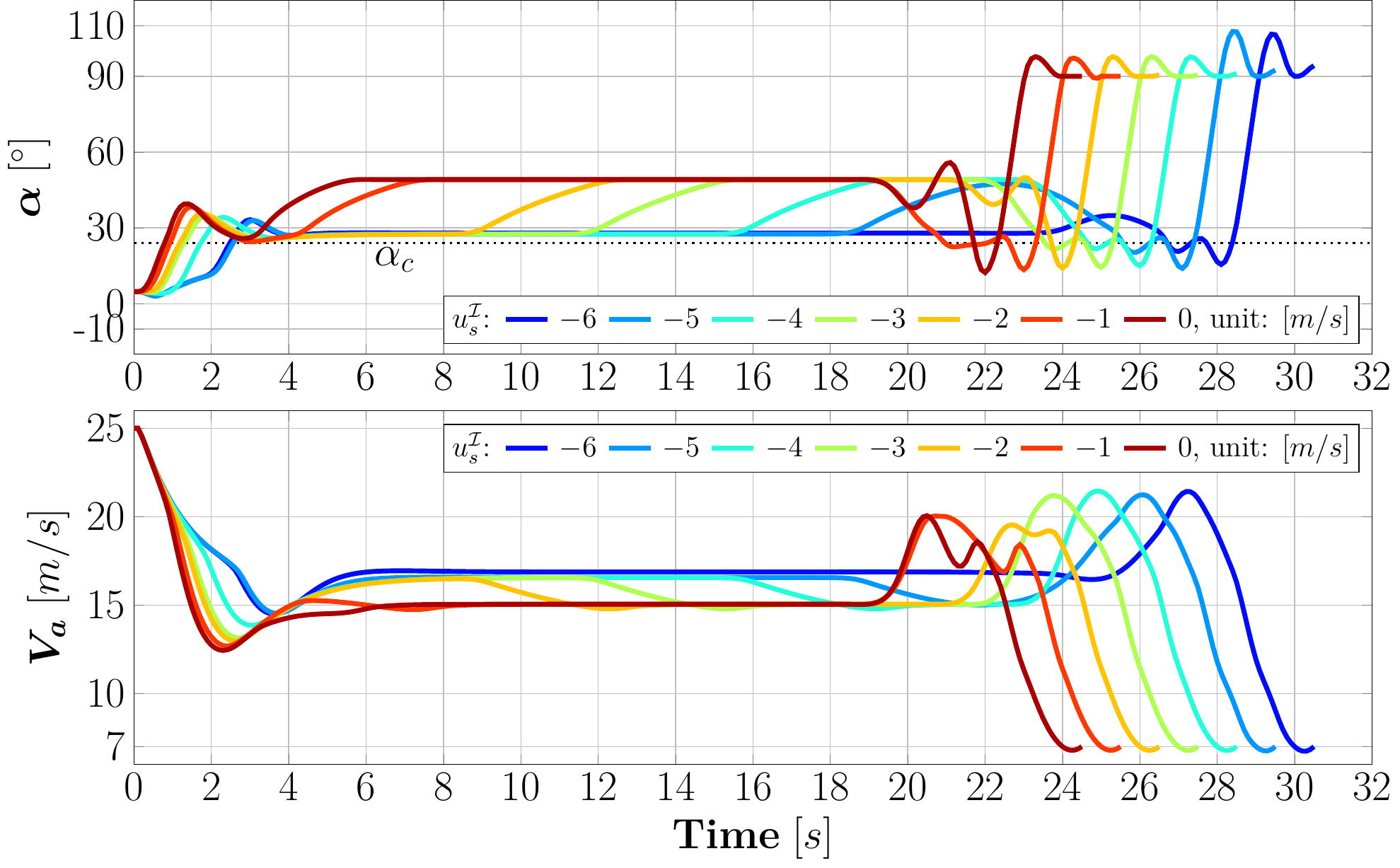}
			\caption{The AoA ($\alpha$) and the pitch ($\theta$) with constant winds.}
			\label{fig_alpha_airspeed_multiple_wind}
		\end{center}
	\end{figure}

	The AoA and the airspeed for various constant wind condition is in Fig. \ref{fig_alpha_airspeed_multiple_wind}. There are two stable deep stall AoAs, with $\alpha = 49.11^\circ$ for the nominal case ($u_s^\mathcal{I}=0m/s$), and $\alpha = 27.93^\circ$ when $u_s^\mathcal{I}=-6m/s$. For $u_s^\mathcal{I} \in [-5,-1] (m/s)$, there is a ``jump'' between those two stable deep stall phases, which leaves an open question for future research. 
\vspace{-15pt}
%%%%%%%%%%%%%%%%%%%%%%%%%%%%%%%%%%%%%%%%%%%%%%%%%%%%%%%%%%%%%
\section{Various initial positions}
\label{sec_appendix_init_pos}
\begin{figure}[htb]
	\begin{center}		 	
		\includegraphics[width=0.8\columnwidth]{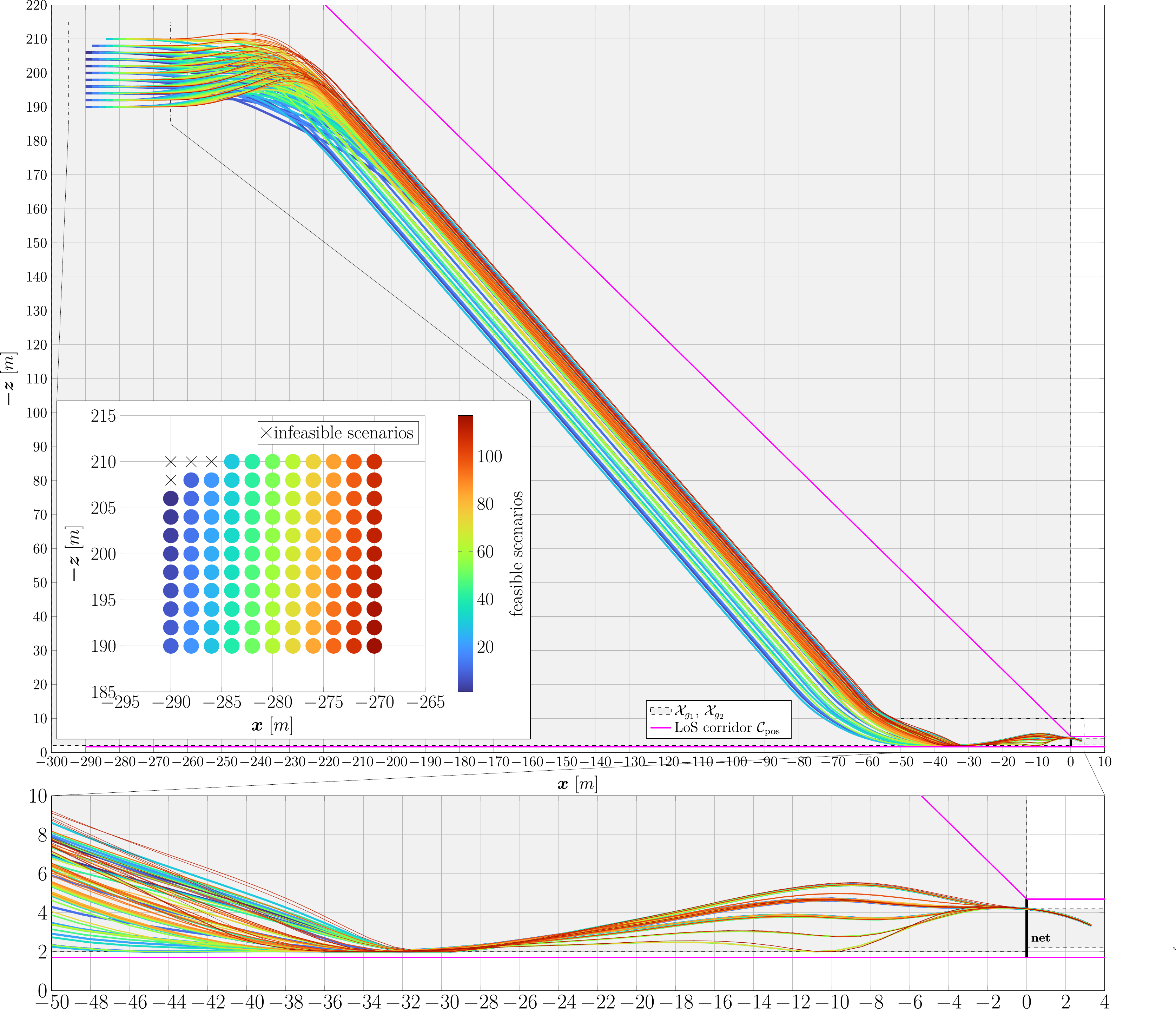}
		\caption{Bundle of trajectory from different initial positions.}
		\label{fig_pos_feasible_infeasible_multiple}
	\end{center}
\end{figure}

\begin{figure}[htb]
	\begin{center}		 	
		\includegraphics[width=0.75\columnwidth]{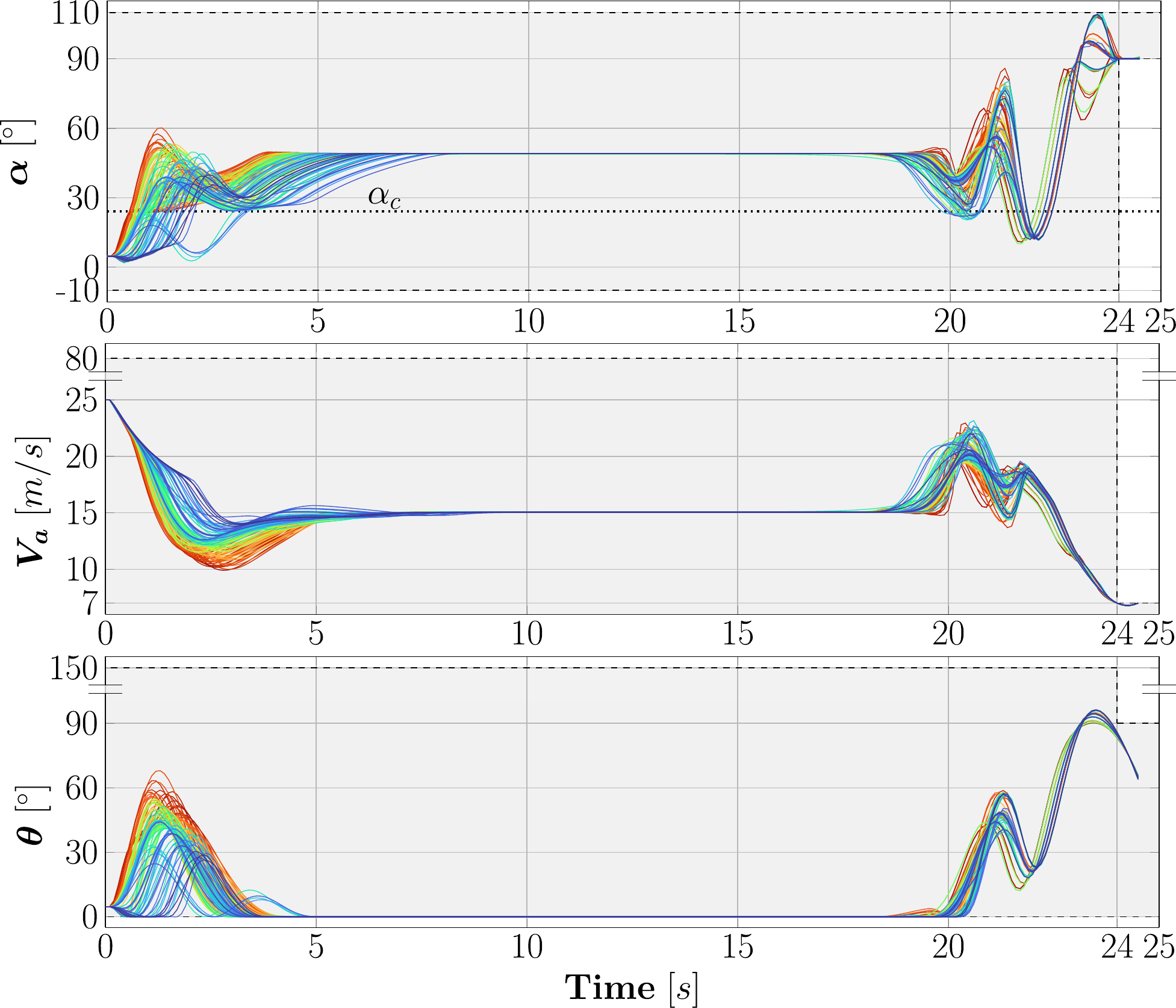}
		\caption{Bundle of AoA, airspeed, and pitch from different initial positions.}
		\label{fig_alpha_airspeed_theta_multiple}
	\end{center}
\end{figure}

\begin{figure}[htb]
	\begin{center}		 	
		\includegraphics[width=0.8\columnwidth]{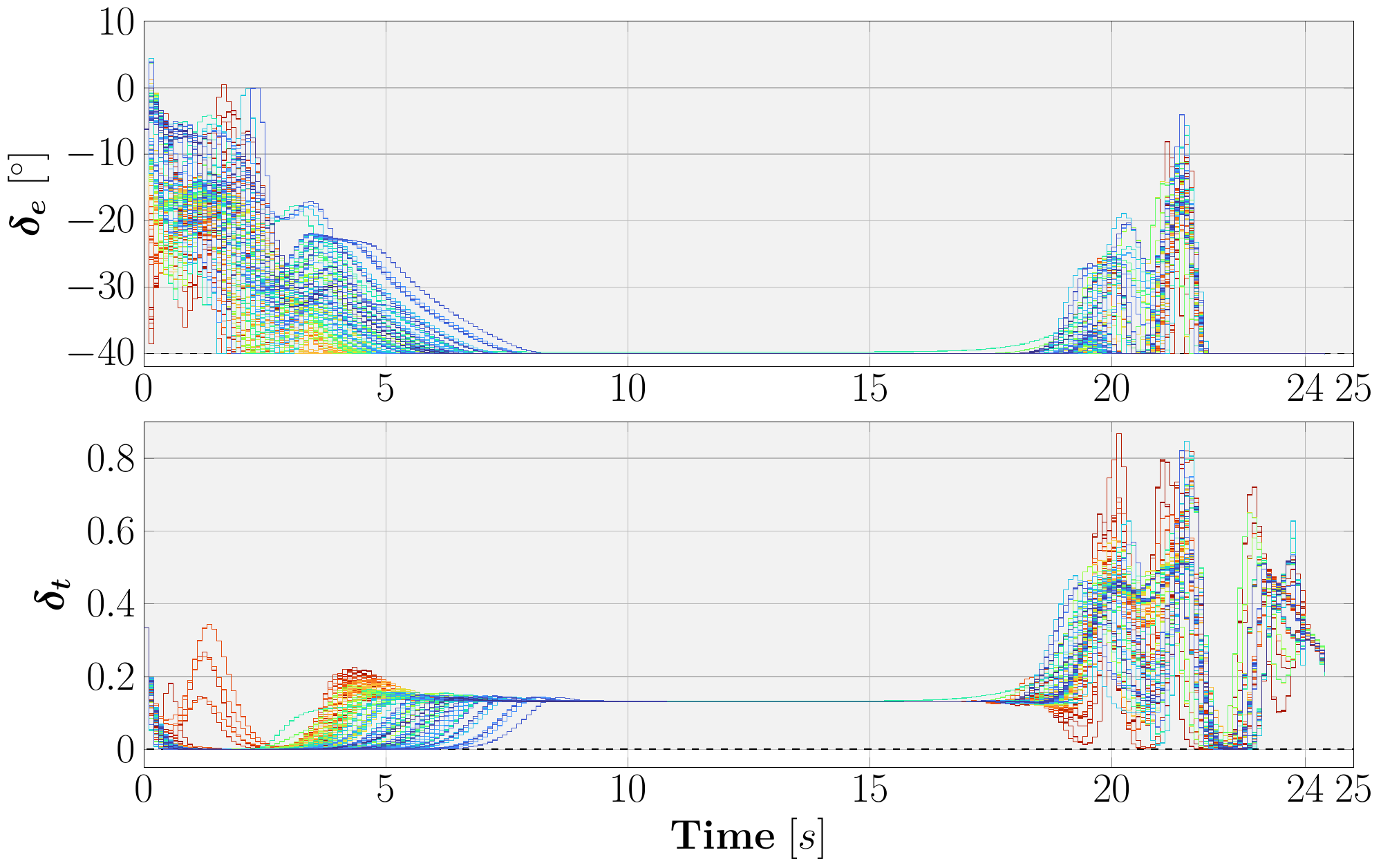}
		\caption{Bundle of control inputs for various initial positions.}
		\label{fig_delta_e_delta_t_multiple}
	\end{center}
\end{figure}

	We solve the generation problem \eqref{eqn_OC_problem_formulation}--\eqref{eqn_OC_problem_formulation_constraints} with the parameters in Subsection \ref{subsec_Simulation_Params}, no wind, by consecutively choose the initial positions in a square $\{x_b,z_b\} \in \{[-290,-270] \times [-210,-190]\} (m)$, with the spacing distance of $2m$, which makes $121$ initial points.   
	\textcolor{black}{The initial guesses (as in Subsection \ref{subsec_Simulation_Params}) for the states are linearly interpolated between 2 points 
		$\bm{\xi}(t_0) = [x_b, \allowbreak z_b,\allowbreak 25,\allowbreak \alpha_0,\allowbreak \theta_0,\allowbreak0,\allowbreak 0,\allowbreak0]^\intercal 
		\text{and} \ 
		\bm{\xi}(t_f + T_p)= [0,\allowbreak-3.2,\allowbreak7,\allowbreak 110^\circ,\allowbreak 90^\circ,\allowbreak0,\allowbreak 0,\allowbreak0]^\intercal$. }
Out of $121$ scenarios, $117$ scenarios give feasible solutions, which are marked by filled colored circles, while other $4$ infeasible solutions are marked as black crosses in Fig. \ref{fig_pos_feasible_infeasible_multiple}. A direct observations could be made: at the furthest upper left points, there are no feasible solutions. 
The bundle of generated feasible optimal trajectories are plotted together in Fig. \ref{fig_pos_feasible_infeasible_multiple}. At the net ($x_r=0m$), all the trajectories \textcolor{black}{tend to} reach the highest point of $\mathcal{X}_{g_2}$ ($-z_r=4.2m$). The bundle of states and inputs are gathered in Figs. \ref{fig_alpha_airspeed_theta_multiple}--\ref{fig_delta_e_delta_t_multiple}. 
	We can observe that the paths have common segments which is a frequent phenomenon in the solution of the OCP (see e.g. \emph{turnpike property} \cite{trelat_TurnpikeProperty_2015}).		
The upper bound constraint on the airspeed for $t<24s$ in Fig. \ref{fig_alpha_airspeed_theta_multiple} is due to $V_a=\sqrt{u^2+w^2}\leq (|u|+|w|) \leq 80 (m/s)$.

\end{appendices}
%%%%%%%%%%%%%%%%%%%%%%%%%%%%%%%%%%%%%%%%%%%%%%%%%%%%%%%%%%%%%
%%%%%%%%%%%%%%%%%%%%%%%%%%%%%%%%%%%%%%%%%%%%%%%%%%%%%%%%%%%%%%
% \addtolength{\textheight}{-12cm}   % This command serves to balance the column lengths
                                  % on the last page of the document manually. It shortens
                                  % the textheight of the last page by a suitable amount.
                                  % This command does not take effect until the next page
                                  % so it should come on the page before the last. Make
                                  % sure that you do not shorten the textheight too much.

%%%%%%%%%%%%%%%%%%%%%%%%%%%%%%%%%%%%%%%%%%%%%%%%%%%%%%%%%%%%%%%%%%%%%%%%%%%%%%%%

%%%%%%%%%%%%%%%%%%%%%%%%%%%%%%%%%%%%%%%%%%%%%%%%%%%%%%%%%%%%%%%%%%%%%%%%%%%%%%%%
\bibliographystyle{IEEEtran}
\bibliography{IEEEabrv,references,ref_deepstall,ref_perching,ref_MPC,ref_ECC_2023}

\end{document}